\documentclass[10pt,a4paper]{article}

\usepackage{amsmath}
\usepackage{graphicx,latexsym,euscript,makeidx,color}
\usepackage{amsmath,amsfonts,amssymb,amsthm,mathrsfs}
\usepackage{enumerate}

\usepackage{amssymb}
\usepackage{amsmath}
\usepackage{amsfonts}
\usepackage{graphicx}
\usepackage{graphicx}          
\usepackage{color}

\usepackage[colorlinks,linkcolor=blue,anchorcolor=green,citecolor=red]{hyperref}

\usepackage{geometry}
\def\5n{\negthinspace \negthinspace \negthinspace \negthinspace \negthinspace }
\def\4n{\negthinspace \negthinspace \negthinspace \negthinspace }
\def\3n{\negthinspace \negthinspace \negthinspace }
\def\2n{\negthinspace \negthinspace }
\def\1n{\negthinspace }

   \def\cF{{\cal F}}

\def\ms{\medskip}

            \def\({\Big (}
                  \def\){\Big )}
          \def\[{\Big[}
           \def\]{\Big]}

\def\bde{\begin{definition}\label}    \def\ede{\end{definition}}
\def\bt{\begin{theorem}\label}        \def\et{\end{theorem}}
\def\bc{\begin{corollary}\label}      \def\ec{\end{corollary}}
\def\bl{\begin{lemma}\label}          \def\el{\end{lemma}}
\def\bp{\begin{proposition}\label}    \def\ep{\end{proposition}}
\def\bas{\begin{assumption}\label}    \def\eas{\end{assumption}}
\def\br{\begin{remark}\label}         \def\er{\end{remark}}
\def\bex{\begin{example}\label}       \def\ex{\end{example}}
\def\ba{\begin{array}}                \def\ea{\end{array}}
\def\be{\begin{equation}}
\def\bel{\begin{equation}\label}      \def\ee{\end{equation}}
\def\bea{\begin{eqnarray*}}           \def\eea{\end{eqnarray*}}

\font\tenbb=msbm10 \font\sevenbb=msbm7 \font\fivebb=msbm5
\newfam\bbfam
\scriptscriptfont\bbfam=\fivebb \textfont\bbfam=\tenbb
\scriptfont\bbfam=\sevenbb

\newtheorem{theorem}{\indent Theorem}[section]
\newtheorem{definition}[theorem]{\indent Definition}
\newtheorem{proposition}[theorem]{\indent Proposition}
\newtheorem{corollary}[theorem]{\indent Corollary}
\newtheorem{lemma}[theorem]{\indent Lemma}
\newtheorem{remark}[theorem]{\indent Remark}
\newtheorem{example}[theorem]{\indent Example}

\newtheorem{assumption}[theorem]{\indent Assumption}

\makeatletter
   
   \@addtoreset{equation}{section}
\makeatother

\allowdisplaybreaks[4]

\def\eqbydef{\mathrel{\stackrel{\Delta}{=}}}

\begin{document}

\title{\bf On Time-Consistent Solution to Time-Inconsistent Linear-Quadratic
Optimal Control of Discrete-Time Stochastic Systems
\thanks{Part of this paper was presented at the 17th IFAC
Symposium on System Identification.}
}
\author{Xun Li\thanks{Department of
Applied Mathematics, The Hong Kong Polytechnic University, Hunghom,
Kowloon, Hong Kong, P.R. China.}~~~~
Yuan-Hua Ni\thanks{College of Computer and Control Engineering,  Nankai University, Tianjin 300350 , P. R. China}~~~~ Ji-Feng Zhang\thanks{Key Laboratory of Systems and
Control, Institute of Systems Science, Academy of Mathematics and
Systems Science, Chinese Academy of Sciences, and the School of Mathematical Sciences, University of Chinese Academy of Sciences,
P. R. China.}}
\maketitle

{\bf Abstract:} In this paper, we investigate a class of time-inconsistent
discrete-time stochastic linear-quadratic optimal control problems,
whose time-consistent solutions consist of an open-loop equilibrium
control and a linear feedback equilibrium strategy. The open-loop
equilibrium control is defined for a given initial pair, while the
linear feedback equilibrium strategy is defined for all the initial
pairs. Maximum-principle-type necessary and sufficient conditions
containing stationary and convexity are derived for the existence of
these two time-consistent solutions, respectively. Furthermore, for
the case where the system matrices are independent of the initial
time, we show that the existence of the open-loop equilibrium
control for a given initial pair is equivalent to the solvability of
a set of nonsymmetric generalized difference Riccati equations and a
set of linear difference equations. Moreover, the existence of
linear feedback equilibrium strategy is equivalent to the
solvability of another set of symmetric generalized difference
Riccati equations.

\ms


\ms


\section{Introduction}

Consider the following discrete-time dynamic system
\begin{eqnarray}\label{system-general}
\left\{\begin{array}{l}X_{k+1}=f(k,X_k,u_k,w_k),\\ [1mm]
X_t=x\in \mathbb{R}^n, ~~k\in\mathbb{T}_t,~~t\in \mathbb{T},
\end{array}
\right.\end{eqnarray}
where $\mathbb{T}=\{0,1,\cdots,N-1\}$ with $N$ being a positive integer, and $\mathbb{T}_t=\{t,\cdots,N-1\}$ for $t\in \mathbb{T}$. In (\ref{system-general}), $\{X_k, k\in \widetilde{{\mathbb{T}}}_t\}$ and $\{u_k, k\in \mathbb{T}_t\}$ with $\widetilde{{\mathbb{T}}}_t=\{t,t+1,\cdots,N\}$ are the state process and the control process, respectively; and $\{w_k,k\in \mathbb{T}_t\}$ is a stochastic disturbance process. Introduce the following cost functional associated with (\ref{system-general})
\begin{eqnarray}\label{performace-general}
J(t,x;u)=\sum_{k=t}^{N-1}\mathbb{E}\big{[}e^{-\delta(k-t)}L(k,X_k,u_k)\big{]}+\mathbb{E}\big{[}e^{-\delta(N-t)}h(X_{N})\big{]}.
\end{eqnarray}
Let $\mathcal{U}[t,N-1]$ be a set of admissible controls. Then, the standard discrete-time stochastic optimal control problem is stated as follows.

\textbf{Problem (C).} For (\ref{system-general}),
(\ref{performace-general}) and the initial pair $(t, x)\in
\mathbb{T}\times \mathbb{R}^n$, find a $\bar{u}\in
\mathcal{U}[t,N-1]$ such that
\begin{eqnarray}\label{eq:Jt}
J(t,x;\bar{u})=\inf_{u\in \mathcal{U}[t,N-1]}J(t,x;u).
\end{eqnarray}
Any $\bar{u}\in\mathcal{U}[t,N-1]$ satisfying (\ref{eq:Jt}) is
called an optimal control for the initial pair $(t, x)$;
$\bar{X}=\{\bar{X}_k=\bar{X}(k;t,x,\bar{u}), k\in
\widetilde{{\mathbb{T}}}_t\}$ is called the corresponding optimal
trajectory, and $(\bar{X},\bar{u})$ is referred as an optimal pair.

For Problem (C), dynamic programming  is a fundamental technique to find its optimal control. 
%
%
By studying the incremental behavior of optimal cost-to-go function as one works backward in time, the derived difference equation is termed as Bellman dynamic programming equation,
%
%
which is associated with Bellman's principle of optimality. Let
$\bar{u}$ be an optimal control of Problem (C) for the initial pair
$(t,x)$, and reconsider versions of Problem (C) along the optimal
trajectory $\bar{X}$. For any $\tau\in
\mathbb{T}_{t+1}=\{t+1,...,N-1\}$, Bellman's principle of optimality
tells us that  $\bar{u}|_{\mathbb{T}_{\tau}}$ (the restriction of
$\bar{u}$ on $\mathbb{T}_{\tau}=\{\tau,...,N-1\}$) is an optimal
control of Problem (C) for the initial pair $(\tau,
\bar{X}(\tau;t,x,\bar{u}))$. This property is also referred as the
time consistency of optimal control. Time consistency ensures that
one needs only to solve an optimal control problem for a given
initial pair, and the obtained optimal control is also optimal along
the optimal trajectory.

However, the time consistency of optimal control fails quite often
in many situations. To see this, let us look at a simple example.


\begin{example}\label{Example-1}  \rm
Consider a stochastic linear-quadratic (LQ, for short) optimal
control with an one-dimensional controlled system
\begin{eqnarray}\label{eq:stoch-ex1}
\left\{\begin{array}{ll}
X_{k+1} = (X_k+u_k) + X_k w_k, \\
X_t = x\in \mathbb{R},~~t=0,~ k \in \{t,...,3\}
\end{array}\right.
\end{eqnarray}
and the cost functional
\begin{eqnarray}\label{eq:stoch-ex2}
J(t,x;u) =\sum_{k=t}^{3}\mathbb{E}\Big{[}\frac{1}{1+(k-t)}u_k^2\Big{]} + \mathbb{E}\Big{[}\frac{2}{1+(4-t)}X_{4}^2\Big{]}.
\end{eqnarray}
Here, $w_0,...,w_{3}\in \mathbb{R}$, are mutually independent with
properties $\mathbb{E}[w_k] = 0$, $\mathbb{E}[w_k^2] = 1$. Find a
control to minimize the cost functional (\ref{eq:stoch-ex2}).

\end{example}

\emph{Solution.} 
Note that $t=0$ and
\begin{eqnarray*}
J(0,x;u)\geq \frac{1}{1+(4-0)}\sum_{k=0}^{3}\mathbb{E}|u_k|^2.
\end{eqnarray*}
Hence, there exists a unique optimal control for the initial pair $(0,x)$,
\begin{eqnarray}\label{Example1.1---1}
\bar{u}^{0,x}_k=-W_{0,k}^{-1}P_{0,k+1}\bar{X}^{0,x}_k\equiv \bar{K}_{0,k}\bar{X}^{0,x}_k,~~k\in \{0,1,2,3\}
\end{eqnarray}
with
\begin{eqnarray*}
\left\{
\begin{array}{l}
P_{0,k}=2P_{0,k+1}-P_{0,k+1}W_{0,k}^{-1}P_{0,k+1},\\
W_{0,k}=\frac{1}{1+(k-0)}+P_{0,k+1},\\
P_{0,4}=\frac{2}{1+(4-0)},~~k\in \{0,1,2,3\}
\end{array}
\right.
\end{eqnarray*}
and
\begin{eqnarray*}
\left\{\begin{array}{ll}
\bar{X}^{0,x}_{k+1} = (\bar{X}^{0,x}_k+\bar{u}^{0,x}_k) + \bar{X}^{0,x}_k w_k, \\
\bar{X}^{0,x}_0 = x,~~ k \in \{0,1,2,3\}.
\end{array}\right.
\end{eqnarray*}
Applying $\bar{u}^{0,x}$, we get $\bar{X}^{0,x}_{1}$. Reconsider the controlled system
\begin{eqnarray*}
\left\{\begin{array}{ll}
X_{k+1} = (X_k+u_k) + X_k w_k, \\
X_{1} = \bar{X}_{1}^{0,x}, ~~ k \in \{1,2,3\}
\end{array}\right.
\end{eqnarray*}
with the cost functional (\ref{eq:stoch-ex2}) (of $t=1$), i.e.,
\begin{eqnarray*}
J(1,\bar{X}^{0,x}_{1};u) =\sum_{k=1}^{3}\mathbb{E}\Big{[}\frac{1}{1+(k-1)}u_k^2\Big{]} + \mathbb{E}\Big{[}\frac{2}{1+(4-1)}X_{4}^2\Big{]}.
\end{eqnarray*}
For such a new LQ problem (with the initial pair $(1, \bar{X}_1^{0,x})$), its optimal control is given by
\begin{eqnarray}\label{Example1.1---2}
\widehat{u}_k=-W_{1,k}^{-1}P_{1,k+1}{\widehat{X}}_k\equiv \widehat{K}_{1,k}{\widehat{X}}_k,~~k\in \{1,2,3\}
\end{eqnarray}
with
\begin{eqnarray*}
\left\{\begin{array}{ll}
\widehat{X}_{k+1} = (\widehat{X}_k+\widehat{u}_k) + \widehat{X}_k w_k, \\
\widehat{X}_{1} = \bar{X}_{1}^{t,x}, ~~ k \in \{1,2,3\}
\end{array}\right.
\end{eqnarray*}
and
\begin{eqnarray*}
\left\{
\begin{array}{l}
P_{1,k}=2P_{1,k+1}-P_{1,k+1}W_{1,k}^{-1}P_{1,k+1},\\
W_{1,k}=\frac{1}{1+(k-1)}+P_{1,k+1},\\
P_{1,N}=\frac{2}{1+(4-1)},~~k\in \{1,2,3\}.
\end{array}
\right.
\end{eqnarray*}
A simple calculation yields
\begin{eqnarray*}
\bar{K}_{0,1}=-0.6038 \mbox{ and } \widehat{K}_{1,1}=-0.4979,
\end{eqnarray*}
which imply $\bar{u}^{0,x}_{1}\neq \widehat{u}_1$.
Hence, the restriction of $\bar{u}^{0,x}$ on $\{1,2,3\}$ is not the
optimal control for the initial pair $(1,\bar{X}_1^{0,x})$ as it is
different from (\ref{Example1.1---2}). Such a phenomenon is referred
as the time inconsistency.
\hfill $\square$

Both $e^{-\delta(k-t)}$ and $\frac{1}{1+(k-t)}$ in the cost
functionals (\ref{performace-general}) and (\ref{eq:stoch-ex2}) are
called the time-discounting functions. The term time discounting is
broadly used ``to encompass any reason for caring less about a
future consequence, including factors that diminish the expected
utility generated by a future consequence" \cite{Frederick}.
$e^{-\delta (k-t)}$ is referred as the constant discounting or
exponential discounting, and $\delta$ is termed as the discounting
rate. As exponential function has the property of group, i.e.,
$e^{-\delta (k-t)}=e^{-\delta (k-\tau)}e^{-\delta (\tau-t)}$,
intertemporal decision optimization problems with exponential
discounting, including Problem (C),  are time-consistent. Such kind
of decision optimization problems is extensively studied in the
communities of economics, finance, and control, etc.

Functions like $\frac{1}{1+(k-t)}$ are of hyperbolic discounting,
which is well documented  and often used to describe the situations
with declining discounting rate. As the hyperbolic discounting
function loses the property of group, the dynamic
optimization problem is time-inconsistent in the sense that
Bellman's principle of optimality no longer holds.
In the survey paper \cite{Frederick}, several types of experimental
evidences are presented to show the reasonableness of hyperbolic
discounting other than exponential discounting; a particular example
is that people always prefer the smaller-sooner reward to the
larger-later reward. Furthermore, quasi-geometric discounting, mean
variance utility and endogenous habit formation are several other
representative examples \cite{Bjork,Krusell} that will ruin the time
consistency.

Introduce the system
\begin{eqnarray}\label{system-general-in}
\left\{\begin{array}{l}
X_{k+1}=f(t,k,x,X_k,u_k,w_k), \\
X_t=x \in \mathbb{R}^n, ~~ k \in\mathbb{T}_t,  ~~ t\in \mathbb{T}
\end{array}\right.
\end{eqnarray}
and the cost functional
\begin{eqnarray}\label{performace-general-in}
J(t,x;u) =\sum_{k=t}^{N-1}\mathbb{E}\big{[}L(t,k,x,X_k,u_k)\big{]}+\mathbb{E}\big{[}h(t,x,X_{N})\big{]}.
\end{eqnarray}
Consider the following optimal control problem.

\textbf{Problem (N).} For  (\ref{system-general-in}),
(\ref{performace-general-in}) and the initial pair  $(t,x)\in
\mathbb{T}\times \mathbb{R}^n$, find a $\bar{u}\in
\mathcal{U}[t,N-1]$ such that
\begin{eqnarray*}
J(t,x;\bar{u})=\inf_{u\in\mathcal{U}[t,N-1]}J(t,x;u),
\end{eqnarray*}
where $\mathcal{U}[t,N-1]$ is a set of admissible controls.

Different from (\ref{system-general}) and
(\ref{performace-general}), the initial time $t$ enters explicitly
into (\ref{system-general-in}) and (\ref{performace-general-in}).
This means the controlled system and the cost functional are
modified at different initial times. Furthermore, the cost
functional (\ref{performace-general-in}) can also be viewed as the
general discounting, which includes the hyperbolic discounting,
exponential discounting and quasi-geometric discounting as special
cases. Therefore, Problem (N) is time-inconsistent, in general.

Due to the time inconsistency, there are two different ways to
handle Problem (N). The first one is the static formulation or
pre-commitment formulation. If the initial strategy is kept all the
way from the initial time, then the strategy can be implemented as
planned. This approach neglects the time inconsistency, and the
optimal control is optimal only when viewed at the initial time.
Differently, another approach addresses the time inconsistency in a
dynamic manner. Instead of seeking an ``optimal control", some kinds
of equilibrium solutions are dealt with. This is mainly motivated by
practical applications in economics and finance, and has recently
attracted considerable interest and efforts.

The explicit formulation of time inconsistency was initiated by
Strotz \cite{Strotz} in 1955, whereas its qualitative analysis can
be traced back to the work of Smith \cite{Smith}. In the discrete-time
case, Strotz's idea is to tackle the time inconsistency by a
lead-follower game with hierarchical structure. Particularly,
controls at different time points were viewed as different selves
(players), and every self integrated the policies of his successor
into his own decision. By a backward procedure, the equilibrium
policy (\emph{if exists}) was obtained. Inspired by Strotz and
intending to tackling practical problems in economics and finance,
lots of works were concerned with time inconsistency of dynamic
systems described by ordinary difference or differential equations;
see, for example,
\cite{Ekeland,Ekland-2,Goldman,Krusell,Laibson,Palacios} and
references therein. Unfortunately, as pointed out by Ekeland
\cite{Ekeland,Ekland-2}, Strotz's equilibrium policy typically fails
to exist and is hard to prove the existence. Thus, it is of great
importance to develop a general theory on time-inconsistent optimal
control. This, on the one hand, can enrich the optimal control
theory, and on the other hand, can provide instructive methodology
to push the solvability of practical problems. Recently, this topic
has attracted considerable attention from the theoretic control
community; see, for example,
\cite{Bjork,Hu-jin-Zhou,Hu-jin-zhou-2,Wang-haiyang,Yong-1,Yong-0,Yong-2013}
and references therein.

Concerned with the time-inconsistent LQ problems, we study two kinds of
time-consistent equilibrium solutions, which are the
open-loop equilibrium control and the closed-loop equilibrium
strategy. The separate investigations of such two formulations are
due to the fact that in the dynamic game theory, open-loop control
distinguishes significantly from closed-loop strategy
\cite{Basar,Bernhard1979,Sun-jingrui,Yong-game-2015}. To compare,
open-loop formulation is to find an open-loop equilibrium
``control", while the ``strategy" is the object of closed-loop
formulation. By a strategy, we mean a decision rule that a
controller uses to select a control action based on the available
information set. Mathematically, a strategy is a mapping or operator
on the information set. When substituting the available information
into a strategy, the open-loop value or open-loop realization of
this strategy is obtained. Furthermore, an open-loop equilibrium
control is corresponding to a given initial pair, whereas the linear
feedback equilibrium strategy  is defined for all the initial pairs.
For more about the terms of open-loop, closed-loop, control and
strategy, we refer the readers to, for example, the monograph
\cite{Basar}.

Strotz's equilibrium solution \cite{Strotz} is essentially a
closed-loop equilibrium strategy, which is further elaborately
developed by Yong to the LQ optimal control \cite{Yong-1,Yong-2013}
as well as the nonlinear optimal control \cite{Yong-0,Yong-2}. In
contrast, open-loop equilibrium control is extensively studied by
Hu-Jin-Zhou \cite{Hu-jin-Zhou,Hu-jin-zhou-2} and Yong
\cite{Yong-2013}. In particular, the closed-loop
formulation can be viewed as the extension of Bellman's dynamic
programming, and the corresponding equilibrium strategy (\emph{if
exists}) is derived by a backward procedure
\cite{Yong-1,Yong-2,Yong-0,Yong-2013}. Differently, the open-loop
equilibrium control is characterized via the maximum-principle-like
methodology \cite{Hu-jin-Zhou,Hu-jin-zhou-2}.

Though the time-inconsistent optimal control has gained considerable
attention, its theory is far from being mature. Concerned with the
LQ problems, general necessary and sufficient conditions still do
not include the existence of time-consistent equilibrium
control/strategy. Hence, more elaborate efforts should be paid on
such topic, and more insightful results are much desirable. In this
paper, a general discrete-time time-inconsistent stochastic LQ
optimal control is investigated, and no definiteness constraint is
posed on the state and control weighting matrices. Such indefinite
setting provides a maximal capacity to model and deal with LQ-type
problems, whose study will generalize existing results to some
extent. For more about standard (time-consistent) indefinite LQ
problems, readers are referred to
\cite{Ait-Chen-Zhou-2002,Ait-Moore-Zhou-1,Chen}. Under this general
condition, this paper intends obtaining some neat results on the
existence of time-consistent equilibrium control/strategy.

The remainder of this paper is organized as follows. In Section 2,
the open-loop time-consistent equilibrium control of Problem (LQ)
for a given initial pair is introduced, whose existence is
characterized by some maximum-principle-type conditions. The
existence of an open-loop equilibrium control for any given initial
pair is then shown to be equivalent to the solvability of a set of
nonsymmetric generalized difference Riccati equations (GDREs, for
short) and a set of linear difference equations (LDEs, for short).
In Section 3, the linear feedback time-consistent equilibrium
strategy is investigated, which is defined for all the initial
pairs.
By also a maximum-principle-like methodology, a set of symmetric
GDREs is introduced to characterized the existence of linear
feedback time-consistent equilibrium strategy.
Sections 4 presents some comparisons between the open-loop
equilibrium control and the linear feedback equilibrium strategy.
Section 5 gives several examples, and Section 6 concludes the paper.

\section{Open-loop time-consistent equilibrium control}\label{Section--open-loop}

Consider the following controlled stochastic difference equation (S$\Delta$E, for short)
\begin{eqnarray}\label{system-1}
\left\{\begin{array}{l}
X^t_{k+1} = A_{t,k}X^t_{k}+B_{t,k}u_k+\big{(}C_{t,k}X^t_{k} +D_{t,k}u_k\big{)}w_k,\\
X^t_{t} = x,~~ k \in  \mathbb{T}_t,~~t\in \mathbb{T},
\end{array}\right.
\end{eqnarray}
where
$A_{t,k},C_{t,k}\in \mathbb{R}^{n\times n}$,
$B_{t,k},D_{t,k}\in \mathbb{R}^{n\times m}$ are deterministic matrices; $\{X^t_{k}, k\in \bar{\mathbb{T}}_t\}\triangleq X^t$ and $\{u_k, k\in \mathbb{T}_t\}\triangleq u$
are the state process and the control process, respectively. In (\ref{system-1}), the initial time $t$ in these matrices and the state is to emphasize the property that the matrices and the state may change according to $t$.
The noise $\{w_k, k\in \mathbb{T}\}$ is assumed to be  a martingale difference sequence defined on
a probability space $(\Omega, \mathcal{F}, P)$ with
\begin{eqnarray}\label{w-moment}
\mathbb{E}[w_{k+1}|\cF_k]=0,~~\mathbb{E}[(w_{k+1})^2|\cF_k]=1,~k\geq 0.
\end{eqnarray}
Here, $\mathbb{E}[\,\cdot\,|\cF_k]$ is the conditional mathematical
expectation with respect to $\cF_k=\sigma\{x_0, w_l, l=0, 1,\cdots,k\}$
and $\mathcal{F}_{-1}$ is understood as $\{\emptyset, \Omega\}$.
The cost functional associated with system (\ref{system-1}) is
\begin{eqnarray}\label{cost-1}
&&J(t,x;u)
 =\sum_{k=t}^{N-1}\mathbb{E}\big{[}(X_k^t)^TQ_{t,k}X^t_k + u_k^TR_{t,k}u_k\big{]}+ \mathbb{E}\big{[}(X_N^t)^TG_tX^t_N\big{]}, 
\end{eqnarray}
where $Q_{t,k}, R_{t,k}, k\in \mathbb{T}_t$ and $G_t$ are
deterministic symmetric matrices  of appropriate dimensions.
Different from \cite{Hu-jin-Zhou} \cite{Yong-1} \cite{Yong-0}
\cite{Yong-2013}, we do not pose any definiteness constraint on the
state and the control weight matrices.
Let $L^2_\mathcal{F}(\mathbb{T}_t; \mathcal{H})$ be a set of
$\mathcal{H}$-valued  processes such that for any its element
$\nu=\{\nu_k, k\in \mathbb{T}_t\}$ $\nu_k$ is
$\mathcal{F}_{k-1}$-measurable and
$\sum_{k=t}^{N-1}\mathbb{E}|\nu_k|^2<\infty$.
In addition, for $k\in \mathbb{T}_t$, $L^2_\mathcal{F}(k;
\mathcal{H})$ is a set  of $\mathcal{H}$-valued random variables
such that any its element $\xi$ is $\mathcal{F}_{k-1}$-measurable
and $\mathbb{E}|\xi|^2<\infty$.
Throughout this paper, $(t,x)$ is called an admissible initial
time-state pair or simply an initial pair for (\ref{system-1}) if
$t\in \mathbb{T}$ and $x\in L^2_{\mathcal{F}}(t;\mathbb{R}^n)$.
Consider the following time-inconsistent stochastic LQ problem.

\textbf{Problem (LQ).}  {For (\ref{system-1}), (\ref{cost-1}) and
the initial pair $(t, x)$, find a ${u}^*\in
L^2_\mathcal{F}(\mathbb{T}_t; \mathbb{R}^m)$, such that
\begin{eqnarray}\label{Problem-LQ}
J(t,x;{u}^*) = \inf_{u\in L^2_\mathcal{F}(\mathbb{T}_t; \mathbb{R}^m)}J(t,x;u).
\end{eqnarray}}

In this section, we investigate the open-loop equilibrium control,
whose definition below is a discrete-time version of
\cite{Hu-jin-Zhou}.

\begin{definition}\label{Definition-open-loop}
$u^{t,x,*}\in L^2_\mathcal{F}(\mathbb{T}_t; \mathbb{R}^m)$ is called
an open-loop equilibrium  control of Problem (LQ) for the initial
pair $(t,x)$, if the following inequality holds for any $k\in
\mathbb{T}_t$ and $u_k\in L^2_\mathcal{F}(k; \mathbb{R}^m)$
\begin{eqnarray}\label{open-loop-equilibrium}
J(k,X^{t,x,*}_k; u^{t,x,*}|_{\mathbb{T}_k})\leq J(k,X^{t,x,*}_k;(u_k, u^{t,x,*}|_{\mathbb{T}_{k+1}})).
\end{eqnarray}
Here, $u^{t,x,*}|_{\mathbb{T}_k}$ and $u^{t,x,*}|_{\mathbb{T}_{k+1}}$ are the restrictions of $u^{t,x,*}$ on $\mathbb{T}_k$ and $\mathbb{T}_{k+1}$, respectively, and $X^{t,x,*}$ is given by
\begin{eqnarray}\label{open-loop-equilibrium-state}
\left\{\begin{array}{l}
X^{t,x,*}_{k+1} = A_{k,k}X^{t,x,*}_{k}+B_{k,k}u^{t,x,*}_k+\big{(}C_{k,k}X^{t,x,*}_{k} +D_{k,k}u^{t,x,*}_k\big{)}w_k,\\[1mm]
X^{t,x,*}_{t} = x,~~ k \in  \mathbb{T}_t.
\end{array}\right.
\end{eqnarray}

\end{definition}

\begin{remark}
By Definition \ref{Definition-open-loop}, an open-loop equilibrium
control $u^{t,x,*}$ is time-consistent along the equilibrium
trajectory $X^{t,x,*}$ in the sense that for any $k\in
\mathbb{T}_t$, $u^{t,x,*}|_{\mathbb{T}_k}$ is an open-loop
equilibrium control for the initial pair $(k,X^{t,x,*}_k)$.
Furthermore, (\ref{open-loop-equilibrium}) is a local optimality condition,
as the control $(u_k, u^{t,x,*}|_{\mathbb{T}_{k+1}})$ differs from $u^{t,x,*}|_{\mathbb{T}_k}=(u^{t,x,*}_k, u^{t,x,*}|_{\mathbb{T}_{k+1}})$ only at time point $k$.
\end{remark}

\begin{remark}
To understand the so-called ``equilibrium" in Definition \ref{Definition-open-loop}, we introduce a game, termed as {Problem (HG)}, with a hierarchical structure. The cost functional of Player $k$ ($k=t,\cdots,N-1$) is
\begin{eqnarray}\label{cost-k}
&&J_k(u_k,u_{-k})
\triangleq J(k,X^k_k; u|_{\mathbb{T}_k})=\sum_{\ell=k}^{N-1}\mathbb{E}\big{[}(X_\ell^k)^TQ_{k,\ell}X_\ell^k + u_{\ell}^TR_{k,\ell}u_{\ell}\big{]}+ \mathbb{E}\big{[}(X_N^k)^TG_{k}X^k_N\big{]} 
\end{eqnarray}
with $u_{-k}\eqbydef\{u_t,\cdots,u_{k-1}, u_{k+1},\cdots, u_{N-1}\}$. In (\ref{cost-k}), $u_k$ is the action of Player $k$, and  $\{X_\ell^k, \ell\in \bar{\mathbb{T}}_k\}\triangleq X^k$ is the internal state of Player $k$ driven by all the actions $\{u_\ell, \ell\in \mathbb{T}_t\}$. Indeed, $\{u_t,\cdots,u_{k-1}\}$ enters into $X^k$ via its initial state
\begin{eqnarray}\label{system-k}
\left\{\begin{array}{l}
X^k_{\ell+1} = A_{k,\ell}X_\ell^k+B_{k,\ell}u_{\ell}+\big{(}C_{k,\ell}X^k_{\ell} +D_{k,\ell}u_{\ell}\big{)}w_{\ell},\\
X_{k}^k = X^{k-1}_k,~~ {\ell} \in  \mathbb{T}_k
\end{array}\right.
\end{eqnarray}
with $X_{k}^{k-1}$ being the internal state of Player $k-1$ at time point $k$, as $X^{k-1}_k$ is essentially a functional of $\{u_t,\cdots, u_{k-1}\}$. This is why we denote $J(k,X_k^k;u|_{\mathbb{T}_k})$ as $J_k(u_k,u_{-k})$.
Furthermore, in (\ref{system-k}), $X_{k}^k = X^{k-1}_k$ indicates a forward hierarchical structure of Problem (HG): Player $k-1$ could be viewed as the leader of Player $k$.
%
%
%
From (\ref{cost-1}), (\ref{open-loop-equilibrium}) and (\ref{cost-k}), we have
\begin{eqnarray*}\label{open-loop-equilibrium-2}
J_k(u_k^{t,x,*}, u_{-k}^{t,x,*})\leq J_k(u_k, u_{-k}^{t,x,*}).
\end{eqnarray*}
Therefore, $u^{t,x,*}$ is the Nash equilibrium of Problem (HG). Hence, in Definition \ref{Definition-open-loop}, we call $u^{t,x,*}$ the equilibrium control.
\end{remark}

The following theorem is concerned with the existence of open-loop equilibrium control.

\begin{theorem}\label{Theorem-Equivalentce-open-loop}
Given an initial pair $(t,x)$, the following statements are equivalent.

(i) There exists an open-loop equilibrium control of Problem (LQ) for the initial pair $(t,x)$.

(ii) There exists a control $u^{t,x,*}$ such that for any $k\in
\mathbb{T}_t$, the following forward-backward stochastic difference
equations (FBS$\Delta$E, for short) has a solution
$(X^{k,*},Z^{k,*})$
\begin{eqnarray}\label{system-adjoint}
\left\{\begin{array}{l}
X^{k,*}_{\ell+1} = A_{k,\ell}X^{k,*}_{\ell}+B_{k,\ell}u^{t,x,*}_\ell +\big{(}C_{k,\ell}X^{k,*}_{\ell} +D_{k,\ell}u^{t,x,*}_\ell\big{)}w_\ell, \\ [1mm]
Z_\ell^{k,*}=A_{k,\ell}^T\mathbb{E}(Z_{\ell+1}^{k,*}|\mathcal{F}_{\ell-1})+C_{k,\ell}^T\mathbb{E}(Z_{\ell+1}^{k,*}w_\ell|\mathcal{F}_{\ell-1})+Q_{k,\ell}X_\ell^{k,*}, \\ [1mm]
X^{k,*}_k=X^{t,x,*}_k,~~Z_N^{k,*}=G_kX^{k,*}_N,~~\ell\in \mathbb{T}_k
\end{array}\right.
\end{eqnarray}
with the stationary condition
\begin{eqnarray}\label{stationary-condition}
0=R_{k,k}u_k^{t,x,*}+B_{k,k}^T\mathbb{E}(Z_{k+1}^{k,*}|\mathcal{F}_{k-1})+D_{k,k}^T\mathbb{E}(Z_{k+1}^{k,*}w_k|\mathcal{F}_{k-1}),
\end{eqnarray}
and the convexity condition
\begin{eqnarray}\label{convex}
\inf_{\bar{u}_k\in L^2_\mathcal{F}(k; \mathbb{R}^m)}\Big{\{}\sum_{\ell=k}^{N-1}\mathbb{E}\big{[}\big{(}Y^k_\ell\big{)}^T Q_{k,\ell} {Y}_\ell^{k}\big{]}+\mathbb{E}\big{[}\bar{u}_k^TR_{k,k}\bar{u}_k\big{]}+\mathbb{E}\big{[}\big{(}Y_N^k\big{)}^T G_{k}{Y}_N^{k}\big{]}\Big{\}}\geq 0.
\end{eqnarray}
Here, $Y^k$ is given by
\begin{eqnarray}\label{system-y-0}
\left\{\begin{array}{l}
{Y}^k_{\ell+1}=A_{k,\ell}{Y}^k_\ell+C_{k,\ell}{Y}^k_\ell w_\ell,~~~\ell\in \mathbb{T}_{k+1},\\[1mm]
Y^k_{k+1}=B_{k,k}\bar{u}_k+D_{k,k}\bar{u}_kw_k,\\[1mm]
{Y}^k_k=0,
\end{array}\right.
\end{eqnarray}
and $X^{t,x,*}$ is given in (\ref{open-loop-equilibrium-state}).

Furthermore, $u^{t,x,*}$ given in (ii) is an open-loop equilibrium control.
\end{theorem}

\emph{Proof}. See Appendix A. \hfill$\square$

To proceed, we now recall the pseudo-inverse of a matrix. By \cite{Penrose}, for a given matrix
$M\in \mathbb{R}^{n\times m}$, there exists a unique matrix in
$\mathbb{R}^{m\times n}$ denoted by $M^\dagger$ such that
\begin{eqnarray}
\left\{
\begin{array}{l}
MM^\dagger M=M,~~ M^\dagger M M^\dagger=M^\dagger,\\
(MM^\dagger)^T=MM^\dagger, ~~(M^\dagger M)^T=M^\dagger M.
\end{array}
\right.
\end{eqnarray}
This $M^\dagger$ is called the Moore-Penrose inverse of $M$. The following lemma is from \cite{Ait-Chen-Zhou-2002}.

\begin{lemma}\label{Lemma-matrix-equation}
Let matrices $L$, $M$ and $N$ be given with appropriate size. Then,
$LXM=N$ has a solution $X$ if and only if $LL^\dagger NMM^\dagger=N$.
Moreover, the solution of $LXM=N$ can be expressed as
$X=L^\dagger NM^\dagger+Y-L^\dagger LYMM^\dagger$,
where $Y$ is a matrix with appropriate size.
\end{lemma}

From (\ref{stationary-condition}) and Lemma \ref{Lemma-matrix-equation}, an open-loop equilibrium control is given by
\begin{eqnarray}\label{optimal-control-1}
u_k^{t,x,*}=-R^{\dagger}_{k,k}\big{[}B_{k,k}^T\mathbb{E}(Z_{k+1}^{k,*}|\mathcal{F}_{k-1})+D_{k,k}^T\mathbb{E}(Z_{k+1}^{k,*}w_k|\mathcal{F}_{k-1})\big{]},~~~k\in \mathbb{T}_t.
\end{eqnarray}
%
%
Substituting (\ref{optimal-control-1}) into (\ref{system-adjoint}), we get a set of FBS$\Delta$Es
\begin{eqnarray}\label{system-adjoint-system}
\left\{\begin{array}{l}
\left\{\begin{array}{l}
X^{k,*}_{\ell+1} =A_{k,\ell}X^{t,*}_{\ell}-B_{k,\ell}R^{\dagger}_{\ell,\ell}\big{[}B_{\ell,\ell}^T\mathbb{E}(Z_{\ell+1}^{\ell,*}|\mathcal{F}_{\ell-1}) +D_{\ell,\ell}^T\mathbb{E}(Z_{\ell+1}^{\ell,*}w_\ell|\mathcal{F}_{\ell-1})\big{]}\\[2mm]
\hphantom{X^{k,*}_{\ell+1} =}+\Big{\{}C_{k,\ell}X^{k,*}_{\ell} -D_{k,\ell}R^{\dagger}_{\ell,\ell}\big{[}B_{\ell,\ell}^T\mathbb{E}(Z_{\ell+1}^{\ell,*}|\mathcal{F}_{\ell-1})+D_{\ell,\ell}^T\mathbb{E}(Z_{\ell+1}^{\ell,*}w_\ell|\mathcal{F}_{\ell-1})\big{]}\Big{\}}w_\ell, \\
[2mm]
Z_\ell^{k,*}=A_{k,\ell}^T\mathbb{E}(Z_{\ell+1}^{k,*}|\mathcal{F}_{\ell-1})+C_{k,\ell}^T\mathbb{E}(Z_{\ell+1}^{k,*}w_\ell|\mathcal{F}_{\ell-1})+Q_{k,\ell}X_\ell^{k,*},\\[2mm]
X^{k,*}_t=X^{t,x,*}_k,~~Z_N^{k,*}=G_kX^{k,*}_N,~~\ell\in \mathbb{T}_{k+1},
\end{array}\right. \\ [2mm]
k\in \mathbb{T}_t,
\end{array}\right.
\end{eqnarray}
coupled with
\begin{eqnarray*}
\left\{\begin{array}{l}
X^{t,x,*}_{k+1} =A_{k,k}X^{t,x,*}_{k}-B_{k,k}R^{\dagger}_{k,k}\big{[}B_{k,k}^T\mathbb{E}(Z_{k+1}^{k,*}|\mathcal{F}_{k-1})+D_{k,k}^T\mathbb{E}(Z_{k+1}^{k,*}w_k|\mathcal{F}_{k-1})\big{]}\\[2mm]
\hphantom{X^{t,x,*}_{k+1} = }+\Big{\{}C_{k,k}X^{t,x,*}_{k} -D_{k,k}R^{\dagger}_{k,k}\big{[}B_{k,k}^T\mathbb{E}(Z_{k+1}^{k,*}|\mathcal{F}_{k-1})+D_{k,k}^T\mathbb{E}(Z_{k+1}^{k,*}w_k|\mathcal{F}_{k-1})\big{]}\Big{\}}w_k,\\[2mm]
X^{t,x,*}_{t} = x,~~ k \in  \mathbb{T}_t.
\end{array}\right.
\end{eqnarray*}
To get a more convenient form of the open-loop equilibrium control, we should decouple (\ref{system-adjoint-system}) to obtain some Riccati-like equations. However, (\ref{system-adjoint-system}) is not a single
FBS$\Delta$E but a set of  FBS$\Delta$Es coupled with the open-loop equilibrium state $X^{t,x,*}$. To this end, we will focus on a specific case where the system equation is
\begin{eqnarray}\label{system-3}
\left\{
\begin{array}{l}
X_{k+1} = (A_{k}X_{k}+B_{k}u_k)+(C_{k}X_{k} +D_{k}u_k)w_k,\\
X_{t} = x,~~ k \in  \mathbb{T}_t.
\end{array}
\right.
\end{eqnarray}
%
%
%
Different from (\ref{system-1}), the system matrices in (\ref{system-3}) are assumed to be independent of the initial time $t$. 
The following result is on the equivalent characterization of the existence of open-loop equilibrium control of Problem (LQ) corresponding to (\ref{system-3}) and (\ref{cost-1}).

\begin{theorem}\label{Theorem-Equivalentce-open-loop-s-g}
The following statements are equivalent.

(i) For any initial pair $(t,x) \in \mathbb{T}\times \mathbb{R}^n$,
there exists an open-loop equilibrium control of Problem (LQ)
corresponding to (\ref{system-3}), (\ref{cost-1}) and the initial
pair $(t,x)$.

(ii) The set of GDREs %
\begin{eqnarray}\label{GDRE-system2}
\left\{\begin{array}{l}
\left\{\begin{array}{l}
P_{t,k}=Q_{t,k}+A_k^T P_{t,k+1}A_k+C_k^T P_{t,k+1}C_k-{H}_{t,k}^T{W}_{k,k}^{\dag}{H}_{k,k},\\[1mm]
P_{t,N}=G_t,\\[1mm]
k\in \mathbb{T}_t,
\end{array}
\right.\\[1mm]
{W}_{t,t}{W}_{t,t}^\dag{H}_{t,t}-{H}_{t,t}=0,\\[1mm]
t\in \mathbb{T}
\end{array}\right.
\end{eqnarray}
and the set of LDEs
\begin{eqnarray}\label{LRE}
\left\{\begin{array}{l}
\left\{\begin{array}{l}
S_{t,k}=Q_{t,k}+A_k^T S_{t,k+1}A_k+C_k^T S_{t,k+1}C_k,\\[1mm]
S_{t,N}=G_t,\\[1mm]
k\in \mathbb{T}_t
\end{array}\right.\\ [1mm]
R_{t,t}+B^T_{t}S_{t,t+1}B_{t}+D^T_{t}S_{t,t+1}D_{t}\geq 0,\\[1mm]
t\in \mathbb{T}
\end{array}\right.
\end{eqnarray}
are solvable in the sense that the constrained conditions ${W}_{t,t}{W}_{t,t}^\dag{H}_{t,t}-{H}_{t,t}=0$, $R_{t,t}+B^T_{t}S_{t,t+1}B_{t}+D^T_{t}S_{t,t+1}D_{t}\geq 0$, $t\in \mathbb{T}$, are satisfied.
Here,
\begin{eqnarray}\label{GDRE-W-H}
\left\{\begin{array}{l}
{W}_{k,k}=R_{k,k}+B^T_{k}P_{k,k+1}B_{k}+D^T_{k}P_{k,k+1}D_{k},\\[1mm]
{H}_{k,k}=B^T_{k}P_{k,k+1}A_{k}+D^T_{k}P_{k,k+1}C_{k},\\[1mm]
H_{t,k}=B^T_{k}P_{t,k+1}A_{k}+D^T_{k}P_{t,k+1}C_{k},\\[1mm]
k\in \mathbb{T}_t, ~t\in \mathbb{T}.
\end{array}\right.
\end{eqnarray}

Under any of above conditions, the open-loop equilibrium control for the initial pair $(t,x)$
admits the feedback form
\begin{eqnarray}\label{Theorem-Equivalentce-open-loop-control}
{u}_k^{t,x,*}=-W_{k,k}^\dagger H_{k,k} {X}_{k}^{t,x,*},~~k\in \mathbb{T}_t,
\end{eqnarray}
where
\begin{eqnarray}\label{open-loop-equilibrium-control-state}
\left\{\begin{array}{l}
{X}^{t,x,*}_{k+1} = \big{(}A_{k}-B_kW_{k,k}^\dagger H_{k,k}\big{)}{X}^{t,x,*}_{k}+\big{(}C_{k}-D_kW_{k,k}^\dagger H_{k,k}\big{)}{X}^{t,x,*}_{k}w_k,\\
{X}^{t,x,*}_{t} = x,~~ k \in  \mathbb{T}_t.
\end{array}\right.
\end{eqnarray}

\end{theorem}

\emph{Proof}. (ii)$\Rightarrow$(i). For any given initial pair
$(t,x)\in \mathbb{T}\times \mathbb{R}^n$, let
\begin{eqnarray}\label{Theorem-Equivalentce-open-loop-GDRE-1}
\widetilde{u}_k^{t,x}=-W_{k,k}^\dagger H_{k,k} \widetilde{X}_{k}^{t,x},~~k\in \mathbb{T}_t,
\end{eqnarray}
where
\begin{eqnarray}\label{open-loop-equilibrium-pair-1}
\left\{\begin{array}{l}
\widetilde{X}^{t,x}_{k+1} = \big{(}A_{k}-B_kW_{k,k}^\dagger H_{k,k}\big{)}\widetilde{X}^{t,x}_{k}+\big{(}C_{k}-D_kW_{k,k}^\dagger H_{k,k}\big{)}\widetilde{X}^{t,x}_{k}w_k,\\
\widetilde{X}^{t,x}_{t} = x,~~ k \in  \mathbb{T}_t.
\end{array}\right.
\end{eqnarray}
Then, the following FBS$\Delta$E
\begin{eqnarray}\label{system-adjoint-system-4}
\left\{\begin{array}{l}
\widetilde{X}^{k}_{\ell+1} = A_{\ell}\widetilde{X}^{k}_{\ell}+B_{\ell}\widetilde{u}_\ell^{t,x}+\big{\{}C_{\ell}\widetilde{X}^{k}_{\ell} +D_{\ell}\widetilde{u}_\ell^{t,x}\big{\}}w_\ell, \\
[1mm]
\widetilde{Z}_\ell^{k}=A_{\ell}^T\mathbb{E}(\widetilde{Z}_{\ell+1}^{k}|\mathcal{F}_{\ell-1}) +C_{\ell}^T\mathbb{E}(\widetilde{Z}_{\ell+1}^{k}w_\ell|\mathcal{F}_{\ell-1})+Q_{k,\ell}\widetilde{X}_k^{k},\\[1mm]
\widetilde{X}^{k}_k=\widetilde{X}_k^{t,x},~~\widetilde{Z}_N^{k}=G_k\widetilde{X}^{k}_N,\\[1mm]
\ell\in \mathbb{T}_k
\end{array}
\right.
\end{eqnarray}
has a solution $(\widetilde{X}^{k}, \widetilde{Z}^{k})$, as the backward state $\widetilde{Z}^{k}$ does not appear in the forward S$\Delta$E.
Comparing the forward S$\Delta$E of (\ref{system-adjoint-system-4}) (by substituting $\widetilde{u}^{t,x}$) with (\ref{open-loop-equilibrium-pair-1}), we have
\begin{eqnarray*}
\widetilde{X}^{k}_\ell=\widetilde{X}^{t,x}_\ell, ~~\ell\in \mathbb{T}_k.
\end{eqnarray*}
To apply Theorem \ref{Theorem-Equivalentce-open-loop}, we validate the stationary condition. Noting $\widetilde{Z}^{k}_N=G_N\widetilde{X}^{k}_N$,  we have
\begin{eqnarray}\label{stationary-condition-ell-4-0}
&&\widetilde{Z}_{N-1}^{k}=\big{(}A^T_{N-1}P_{k,N}A_{N-1}+C^T_{N-1}P_{k,N}C_{N-1}+Q_{k,N-1}\big{)}\widetilde{X}_{N-1}^{k}\nonumber\\
&&\hphantom{\widetilde{Z}_{N-1}^{k}=}+\big{(}A^T_{N-1}P_{k,N}B_{N-1}+C^T_{N-1}P_{k,N}D_{N-1} \big{)}\widetilde{u}_{N-1}^{t,x}\nonumber\\
&&\hphantom{\widetilde{Z}^{k}_{N-1}}=\Big{(}A^T_{N-1}P_{k,N}A_{N-1}+C^T_{N-1}P_{k,N}C_{N-1}+Q_{k,N-1}-H_{k,N-1}^TW_{N-1,N-1}^\dagger H_{N-1,N-1}\Big{)}\widetilde{X}^{k}_{N-1}\nonumber\\
&&\hphantom{\widetilde{Z}^{k}_{N-1}}=P_{k,N-1}\widetilde{X}^{k}_{N-1}.
\end{eqnarray}
By some backward induction, we get
\begin{eqnarray*}
\widetilde{Z}^{k}_\ell=P_{k,\ell}\widetilde{X}^{k}_\ell,~~\ell\in \mathbb{T}_k,~~k\in \mathbb{T}_t,
\end{eqnarray*}
which implies
\begin{eqnarray*}
&&R_{k,k}\widetilde{u}_k^{t,x}+B_{k}^T\mathbb{E}(\widetilde{Z}_{k+1}^{k}|\mathcal{F}_{k-1})+D_{k}^T\mathbb{E}(\widetilde{Z}_{k+1}^{k}w_k|\mathcal{F}_{k-1})\\[1mm]
&&=\big{(}R_{k,k}+B^T_{k}P_{k,k+1}B_{k}+D^T_{k}P_{k,k+1}D_{k}\big{)}\widetilde{u}_k^{{t,x}}+\big{(}B^T_{k}P_{k,k+1}A_{k}+D^T_{k}P_{k,k+1}C_{k}\big{)}\widetilde{X}_{k}^{{t,x}}\nonumber\\[1mm]
&&= W_{k,k}\widetilde{u}_k^{t,x}+H_{k,k}\widetilde{X}^{t,x}_k\\[1mm]
&&=0.
\end{eqnarray*}
The last equality follows from the solvability of
(\ref{GDRE-system2}).
Therefore, the stationary condition holds.
Furthermore, corresponding to (\ref{system-y-0}) and (\ref{system-3}), let
\begin{eqnarray*}
\left\{\begin{array}{l}
{Y}_{\ell+1}=A_{\ell}{Y}_\ell+C_{\ell}{Y}_\ell w_\ell,~~~\ell\in \mathbb{T}_{k+1},\\
Y_{k+1}=B_{k}\bar{u}_k+D_{k}\bar{u}_kw_k,\\
{Y}_k=0.
\end{array}\right.
\end{eqnarray*}
Then, by some simple calculations, we have for any $\bar{u}_k\in L^2_{\mathcal{F}}(k; \mathbb{R}^m)$
\begin{eqnarray}\label{convex-3}
&&\mathbb{E}\big{[} \bar{u}^T_k R_{k,k} \bar{u}_k\big{]}+\mathbb{E}\big{[} (Y)^T_N G_{k}Y_N\big{]}+\sum_{\ell=k}^{N-1}\mathbb{E}\big{[} (Y_\ell)^T Q_{k,\ell}Y_\ell\big{]}\nonumber \\[1mm]
&&=\sum_{\ell=k}^{N-1}\mathbb{E}\Big{[} (Y_{\ell+1})^T S_{k,\ell+1}Y_{\ell+1}-(Y_{\ell})^T (S_{k,\ell}- Q_{k,\ell})Y_\ell\Big{]}+\mathbb{E}\big{[} \bar{u}^T_k R_{k,k} \bar{u}_k\big{]}\nonumber \\[1mm]
&&=\mathbb{E}\big{[}\bar{u}_k^T\big{(}R_{k,k}+B^T_{k}S_{k,k+1}B_{k}+D^T_{k}S_{k,k+1}D_{k}\big{)}\bar{u}_k\big{]}\nonumber\\[1mm]
&&\geq 0,
\end{eqnarray}
where the inequality is from the solvability of (\ref{LRE}). Therefore, the convexity condition holds. By Theorem \ref{Theorem-Equivalentce-open-loop}, the pair $(\widetilde{X}^{t,x}, \widetilde{u}^{t,x})$ given in (\ref{Theorem-Equivalentce-open-loop-GDRE-1})-(\ref{open-loop-equilibrium-pair-1}) is an  open-loop equilibrium pair of Problem (LQ) for (\ref{system-3}), (\ref{cost-1}) and $(t,x)$.

(i)$\Rightarrow$(ii). For any given initial pair $(t,x)\in
\mathbb{T}\times \mathbb{R}^n$, let $(X^{t,x,*}, u^{t,x,*})$ be an
open-loop equilibrium pair of Problem (LQ) corresponding to
(\ref{system-3}) and (\ref{cost-1}). In this case, for any $k\in
\mathbb{T}_t$, (\ref{system-adjoint})  (\ref{stationary-condition})
become
\begin{eqnarray}\label{system-adjoint-system-4-2}
\left\{
\begin{array}{l}
{X}^{t,x,*}_{\ell+1} = A_{\ell}{X}^{t,x,*}_{\ell}+B_{\ell}{u}_\ell^{t,x,*}+\big{(}C_{\ell}{X}^{t,x,*}_{\ell} +D_{\ell}{u}_\ell^{t,x,*}\big{)}w_\ell, \\
[1mm]
{Z}_\ell^{k,*}=Q_{k,\ell}{X}_k^{t,x,*}+A_{\ell}^T\mathbb{E}({Z}_{\ell+1}^{k,*}|\mathcal{F}_{\ell-1}) +C_{\ell}^T\mathbb{E}({Z}_{\ell+1}^{k,*}w_\ell|\mathcal{F}_{\ell-1}),\\[1mm]
{X}^{t,x,*}_k={X}_k^{t,x,*},~~{Z}_N^{k,*}=G_k{X}^{t,x,*}_N,\\[1mm]
\ell\in \mathbb{T}_k,
\end{array}
\right.
\end{eqnarray}
and
\begin{eqnarray}\label{stationary-condition-2}
0=R_{k,k}u_k^{t,x,*}+B_{k}^T\mathbb{E}(Z_{k+1}^{k,*}|\mathcal{F}_{k-1})+D_{k}^T\mathbb{E}(Z_{k+1}^{k,*}w_k|\mathcal{F}_{k-1}).
\end{eqnarray}
Letting $t=N-1$, we have from (\ref{system-adjoint-system-4-2}) and (\ref{stationary-condition-2})
\begin{eqnarray}\label{stationary-condition--open-0}
&&0=R_{N-1,N-1}u_{N-1}^{N-1,x,*}+B_{N-1}^T\mathbb{E}(Z_{N}^{N-1,*}|\mathcal{F}_{N-2})+D_{N-1}^T\mathbb{E}(Z_{N}^{N-1,*}w_{N-1}|\mathcal{F}_{N-2})\nonumber \\[1mm]
&&\hphantom{0}=\big{(}R_{N-1,N-1}+B^T_{N-1}G_{N-1}B_{N-1}+D^T_{N-1}G_{N-1}D_{N-1} \big{)}u_{N-1}^{{N-1,x,*}}\nonumber\\[1mm]
&&\hphantom{0=}\nonumber +\big{(}B^T_{N-1}G_{N-1}A_{N-1}+D^T_{N-1}G_{N-1}C_{N-1}\big{)}X_{N-1}^{{N-1,x,*}}\nonumber \\[1mm]
&&\hphantom{0}\triangleq W_{N-1,N-1}u_{N-1}^{{N-1,x,*}}+H_{N-1,N-1}X^{N-1,x,*}_{N-1},
\end{eqnarray}
where $W_{N-1,N-1}, H_{N-1,N-1}$ are defined in (\ref{GDRE-W-H}) with $P_{N-1,N}=G_{N-1}$. As $x$ in (\ref{stationary-condition--open-0}) can be taken arbitrarily, we have from Lemma \ref{Lemma-matrix-equation}
\begin{eqnarray}\label{W-H-1}
{W}_{N-1,N-1}{W}_{N-1,N-1}^\dag{H}_{N-1,N-1}-{H}_{N-1,N-1}=0,
\end{eqnarray}
and
\begin{eqnarray*}
u^{N-1,x,*}_{N-1}=-{W}_{N-1,N-1}^\dag{H}_{N-1,N-1} x.
\end{eqnarray*}

Assume that we have derived the GDREs over the time period $\mathbb{T}_{t+1}=\{t+1,...,N-1\}$, namely,
\begin{eqnarray}\label{GDRE-system2---t+1}
\left\{\begin{array}{l}
\left\{\begin{array}{l}
P_{k,\ell}=Q_{k,\ell}+A_\ell^T P_{k,\ell+1}A_\ell+C_\ell^T P_{k,\ell+1}C_\ell-{H}_{k,\ell}^T{W}_{\ell,\ell}^{\dag}{H}_{\ell,\ell},\\[1mm]
P_{k,N}=G_k,\\[1mm]
\ell\in \mathbb{T}_k,
\end{array}
\right.\\[1mm]
{W}_{k,k}{W}_{k,k}^\dag{H}_{k,k}-{H}_{k,k}=0,\\[1mm]
k\in \mathbb{T}_{t+1}
\end{array}\right.
\end{eqnarray}
is solvable. To prove the solvability of the GDREs over $\mathbb{T}_{t}$, it is necessary to prove that the GDRE associated with $t$
\begin{eqnarray}\label{GDRE-system2---t}
\left\{\begin{array}{l}
\left\{\begin{array}{l}
P_{t,k}=Q_{t,k}+A_k^T P_{t,k+1}A_k+C_k^T P_{t,k+1}C_k-{H}_{t,k}^T{W}_{k,k}^{\dag}{H}_{k,k},\\[1mm]
P_{t,N}=G_t,\\[1mm]
k \in \mathbb{T}_t,
\end{array}
\right.\\[1mm]
{W}_{t,t}{W}_{t,t}^\dag{H}_{t,t}-{H}_{t,t}=0
\end{array}\right.
\end{eqnarray}
is solvable.
Let $x\in \mathbb{R}^n$ and consider Problem (LQ) corresponding to (\ref{system-3}), (\ref{cost-1}) and the initial pair $(t,x)$.
Similarly to (\ref{stationary-condition--open-0}), we have for $k=N-1$
\begin{eqnarray}\label{stationary-condition--open-1}
0=W_{N-1,N-1}u_{N-1}^{{t,x,*}}+H_{N-1,N-1}X^{t,x,*}_{N-1},
\end{eqnarray}
which together with the solvability of (\ref{GDRE-system2---t+1})
implies
\begin{eqnarray}\label{open-loop-equilibrium-control-state-N-1}
u^{t,x,*}_{N-1}=-{W}_{N-1,N-1}^\dag{H}_{N-1,N-1} X^{t,x,*}_{N-1}.
\end{eqnarray}
Similar to (\ref{stationary-condition-ell-4-0}), we have
\begin{eqnarray}\label{Z-N-1}
&&Z^{N-2,*}_{N-1}=\Big{(}A^T_{N-1}G_{N-2}A_{N-1}+C^T_{N-1}G_{N-2}C_{N-1}+Q_{N-2,N-1}\nonumber \\[1mm]
&&\hphantom{Z^{N-2,*}_{N-1}=}-H_{N-2,N-1}^TW_{N-1,N-1}^\dagger H_{N-1,N-1}\Big{)}X^{t,x,*}_{N-1}\nonumber\\[1mm]
&&\hphantom{Z^{N-2,*}_{N-1}}= P_{N-2,N-1}X^{t,x,*}_{N-1}.
\end{eqnarray}
%
%
%
Furthermore, from (\ref{system-adjoint-system-4-2}) and (\ref{stationary-condition-2}) we have
\begin{eqnarray*}\label{stationary-condition--open-2}
&&0=R_{N-2,N-2}u_{N-2}^{t,x,*}+B_{N-2}^T\mathbb{E}(Z_{N-1}^{N-2,*}|\mathcal{F}_{N-3})+D_{N-2}^T\mathbb{E}(Z_{N-1}^{N-2,*}w_{N-2}|\mathcal{F}_{N-3})\nonumber \\[1mm]
&&\hphantom{0}=\big{(}R_{N-2,N-2}+B^T_{N-2}P_{N-2,N-1}B_{N-2}+D^T_{N-2}P_{N-2,N-1}D_{N-2} \big{)}u_{N-2}^{{t,x,*}}\nonumber\\[1mm]
&&\hphantom{0=}\nonumber +\big{(}B^T_{N-2}P_{N-2,N-1}A_{N-2}+D^T_{N-2}P_{N-2,N-1}C_{N-2}\big{)}X_{N-2}^{{t,x,*}}\nonumber \\[1mm]
&&\hphantom{0}= W_{N-2,N-2}u_{N-2}^{{t,x,*}}+H_{N-2,N-2}X^{t,x,*}_{N-2}.
\end{eqnarray*}
By the solvability of (\ref{GDRE-system2---t+1}) and Lemma \ref{Lemma-matrix-equation}, we have
\begin{eqnarray}\label{open-loop-equilibrium-control-state-N-2}
u^{t,x,*}_{N-2}=-{W}_{N-2,N-2}^\dag{H}_{N-2,N-2} X^{t,x,*}_{N-2}.
\end{eqnarray}
Furthermore,
\begin{eqnarray}\label{Z-N-1-(2)}
&&Z^{N-3,*}_{N-1}=\Big{\{}Q_{N-3,N-1}+A^T_{N-1}G_{N-3}A_{N-1}+C^T_{N-1}G_{N-3}C_{N-1}\nonumber \\[1mm]
&&\hphantom{Z^{N-3,*}_{N-1}=}-H_{N-3,N-1}^TW_{N-1,N-1}^\dagger H_{N-1,N-1}\Big{\}}X^{t,x,*}_{N-1}\nonumber\\[1mm]
&&\hphantom{Z^{N-2,*}_{N-1}}= P_{N-3,N-1}X^{t,x,*}_{N-1},
\end{eqnarray}
and
\begin{eqnarray}\label{Z-N-2}
&&Z^{N-3,*}_{N-2}=Q_{N-3,N-2}X^{t,x,*}_{N-2}+A_{N-2}^T\mathbb{E}\big{(}Z^{N-3,*}_{N-1}|\mathcal{F}_{N-3}\big{)}+C_{N-2}^T \mathbb{E}\big{(}Z^{N-3,*}_{N-1}w_{N-2}|\mathcal{F}_{N-3}\big{)}\nonumber \\
&&\hphantom{Z^{N-3,*}_{N-2}}=\big{(}Q_{N-3,N-2}+A^T_{N-2}P_{N-3,N-1}A_{N-2}+C^T_{N-2}P_{N-3,N-1}C_{N-1}\big{)}X^{t,x,*}_{N-2}\nonumber \\[1mm]
&&\hphantom{Z^{N-3,*}_{N-2}=}+\big{(}A_{N-2}^TP_{N-3,N-1}B_{N-2}+C_{N-2}^TP_{N-3,N-1}D_{N-2} \big{)}u^{t,x,*}_{N-2}\nonumber\\[1mm]
&&\hphantom{Z^{N-3,*}_{N-2}}=\Big{(}Q_{N-3,N-2}+A^T_{N-2}P_{N-3,N-1}A_{N-2}+C^T_{N-2}P_{N-3,N-1}C_{N-1}\nonumber \\[1mm]
&&\hphantom{Z^{N-3,*}_{N-2}=}-H_{N-3,N-2}^TW_{N-2,N-2}^\dagger H_{N-2,N-2} \Big{)}X^{t,x,*}_{N-2}\nonumber \\[1mm]
&&\hphantom{Z^{N-2,*}_{N-1}}= P_{N-3,N-2}X^{t,x,*}_{N-2}.
\end{eqnarray}
From (\ref{open-loop-equilibrium-control-state-N-1})-(\ref{Z-N-2}), we have the following deduction
\begin{eqnarray}\label{u-z}
\left\{
\begin{array}{l}
\left\{
\begin{array}{l}
u^{t,x,*}_{\ell}=-{W}_{\ell,\ell}^\dag{H}_{\ell,\ell} X^{t,x,*}_{\ell},\\
Z^{k-1,*}_{\ell}=P_{k-1,\ell}X^{t,x,*}_{\ell},\\
\ell\in \mathbb{T}_{k},
\end{array}
\right.\\
k\in \mathbb{T}_{t+1}.
\end{array}
\right.
\end{eqnarray}
To extend (\ref{u-z}) to the case of including $k=t$, we have from (\ref{system-adjoint-system-4-2}) and (\ref{stationary-condition-2})
\begin{eqnarray*}\label{stationary-condition--open-2}
&&0=R_{t,t}u_{t}^{t,x,*}+B_{t}^T\mathbb{E}(Z_{t+1}^{t,*}|\mathcal{F}_{t-1})+D_{t}^T\mathbb{E}(Z_{t+1}^{t,*}w_{t}|\mathcal{F}_{t-1})\nonumber \\[1mm]
&&\hphantom{0}=\big{(}R_{t,t}+B^T_{t}P_{t,t+1}B_{t}+D^T_{t}P_{t,t+1}D_{t} \big{)}u_{t}^{{t,x,*}}+\big{(}B^T_{t}P_{t,t+1}A_{t}+D^T_{t}P_{t,t+1}C_{t}\big{)}x\nonumber \\[1mm]
&&\hphantom{0}= W_{t,t}u_{t}^{{t,x,*}}+H_{t,t}x.
\end{eqnarray*}
As $x$ can be arbitrarily taken, from Lemma
\ref{Lemma-matrix-equation} we have
\begin{eqnarray}\label{W-H-t}
{W}_{t,t}{W}_{t,t}^\dag{H}_{t,t}-{H}_{t,t}=0,
\end{eqnarray}
and
\begin{eqnarray*}
u^{t,x,*}_{t}=-{W}_{t,t}^\dag{H}_{t,t} x.
\end{eqnarray*}
By (\ref{W-H-t}), the GDRE (\ref{GDRE-system2---t}) associated with
$t$ is solvable. Thus, GDREs over $\mathbb{T}_{t}$ are solvable.
By the method of induction, we have the solvability of the set of GDREs (\ref{GDRE-system2}).

To conclude the proof, we need to show the solvability of (\ref{LRE}). From (\ref{convex}) and similar to (\ref{convex-3}), we have
\begin{eqnarray*}
0\leq \inf_{\bar{u}_k\in L^2_{\mathcal{F}}(k;\mathbb{R}^m)}\mathbb{E}\big{[}\bar{u}_k^T\big{(}R_{k,k}+B^T_{k}S_{k,k+1}B_{k}+D^T_{k}S_{k,k+1}D_{k}\big{)}\bar{u}_k\big{]}, \quad k\in \mathbb{T}.
\end{eqnarray*}
Hence,
\begin{eqnarray*}
R_{k,k}+B^T_{k}S_{k,k+1}B_{k}+D^T_{k}S_{k,k+1}D_{k} \geq 0, ~~~k\in \mathbb{T},
\end{eqnarray*}
and (\ref{LRE}) is solvable.
\hfill $\square$

\begin{remark}
In (\ref{GDRE-system2}), the GDREs are coupled via  $\{W_{k,k}^\dagger H_{k,k},k\in \mathbb{T}\}$. As for $k\neq t$ $H_{k,k}$ is generally not equal to $H_{t,k}$, $P_{t,k}, k\in \mathbb{T}_t, t\in \mathbb{T}$, are nonsymmetric.
On the contrary, the LDEs in (\ref{LRE}) are decoupled.
Hence, $S_{t,k}, k\in \mathbb{T}_t, t\in \mathbb{T}$, are all symmetric as $G_t, t\in \mathbb{T}$, are symmetric.
Interestingly, there is no definite constraint on matrices associated with (\ref{GDRE-system2}), while the definite constraint is posed through (\ref{LRE}). This is different from the standard indefinite stochastic LQ optimal control. Concerned with the reasons, the definite constraints in (\ref{LRE}) are equivalent to the convexity conditions, while the constraints in (\ref{GDRE-system2}) are associated with the stationary conditions.
\end{remark}

Let us further assume that $G_t$ and $Q_{t,k} (k\in \mathbb{T}_t)$ are all independent of $t$.
Then, (\ref{GDRE-system2}) and (\ref{LRE}) become
\begin{eqnarray}\label{GDRE-system2-s}
\left\{\begin{array}{l}
P_{k}=Q_{k}+A_k^T P_{k+1}A_k+C_k^T P_{k+1}C_k -{H}_{k}^T{W}_{k}^{\dag}{H}_{k}, \\ [1mm]
P_{N}=G, \\ [1mm]
{W}_{k}{W}_{k}^\dag{H}_{k}-{H}_{k}=0, \\ [1mm]
k\in \mathbb{T},
\end{array}\right.
\end{eqnarray}
and
\begin{eqnarray}\label{LRE-s}
\left\{\begin{array}{l}
S_{k}=Q_{k}+A_k^T S_{k+1}A_k+C_k^T S_{k+1}C_k,\\[1mm]
S_{N}=G,\\[1mm]
R_{k,k}+B^T_{k}S_{k+1}B_{k}+D^T_{k}S_{k+1}D_{k}\geq 0,\\[1mm]
k\in \mathbb{T},
\end{array}\right.
\end{eqnarray}
where
\begin{eqnarray*}\label{GDRE-W-H-s}
\left\{\begin{array}{l}
{W}_{k}=R_{k,k}+B^T_{k}P_{k+1}B_{k}+D^T_{k}P_{k+1}D_{k},\\[1mm]
{H}_{k}=B^T_{k}P_{k+1}A_{k}+D^T_{k}P_{k+1}C_{k},\\[1mm]
k\in \mathbb{T}.
\end{array}\right.
\end{eqnarray*}
In this case, $P_k, k\in \mathbb{T}$, are all symmetric. Review the standard GDRE \cite{Ait-Chen-Zhou-2002}
\begin{eqnarray}\label{GDRE-s}
\left\{\begin{array}{l}
P_{k}=Q_{k}+A_k^T P_{k+1}A_k+C_k^T P_{k+1}C_k-{H}_{k}^T{W}_{k}^{\dag}{H}_{k}, \\ [1mm]
P_{N}=G,\\[1mm]
{W}_{k}{W}_{k}^\dag{H}_{k}-{H}_{k}=0,\\[1mm]
W_k\geq 0,\\[1mm]
k\in \mathbb{T}.
\end{array}\right.
\end{eqnarray}
Though the GDRE (\ref{GDRE-system2-s}) and the LDE (\ref{LRE-s}) are different from the GDRE (\ref{GDRE-s}), we claim that (\ref{GDRE-system2-s}) and (\ref{LRE-s}) are both solvable for the case with $Q_{k}\geq 0, G\geq 0$ and $R_{k,k}\geq 0$, $k\in \mathbb{T}$.
However, the condition that ensure the solvability of (\ref{GDRE-system2}) is hard to obtain (even for the definite case) due to its nonsymmetric structure. At the present time, we therefore need to validate the solvability of (\ref{GDRE-system2}) case by case.
In the future, we shall study the condition that ensure the solvability of (\ref{GDRE-system2}), and focus on more general cases other than the case of (\ref{system-3}).

\section{Linear feedback time-consistent equilibrium strategy}\label{Section--closed-loop}

In this section, we investigate a kind of closed-loop equilibrium solution, which focuses on the time-consistency of the strategy.
Here, a strategy means a decision rule that a controller uses to select her control action, based on available information set. Mathematically, a strategy is a measurable mapping on the information set. When we substitute the available information into a strategy, the open-loop realization or open-loop value of this strategy is obtained.

\begin{definition}\label{definition-closed-loop}
$\Phi=\{\Phi_0,\cdots, \Phi_{N-1}\}$ with $\Phi_t\in
\mathbb{R}^{m\times n}, t\in \mathbb{T}$, is called a linear
feedback equilibrium strategy of Problem (LQ), if the following
inequality holds for any initial pair $(t,x)\in \mathbb{T} \times
\mathbb{R}^n$, $k\in \mathbb{T}_t$ and $u_k\in L^2_\mathcal{F}(k;
\mathbb{R}^m)$
\begin{eqnarray}\label{defi-closed-loop}
J\big{(}k, X_k^{t,x,*}; (\Phi X^{k,\Phi})|_{\mathbb{T}_k}\big{)}\leq J\big{(}k, X_k^{t,x,*}; (u_k,(\Phi X^{k,u_k,\Phi})|_{\mathbb{T}_{k+1}})\big{)}.
\end{eqnarray}
Here, $X^{t,x,*}=\{X^{t,x,*}_{k}, k\in {\widetilde{\mathbb{T}}}_t\}$, $X^{k,\Phi}=\{X^{k,\Phi}_{\ell}, \ell\in {\widetilde{\mathbb{T}}}_k\}$ and $X^{k,u_k,\Phi}=\{X^{k,u_k,\Phi}_{\ell}, \ell\in {\widetilde{\mathbb{T}}}_k\}$ are given, respectively, by
\begin{eqnarray}
&&\label{equilibrium-state-closed}\left\{\begin{array}{l}
X^{t,x,*}_{k+1} =\big{(}A_{k,k}+B_{k,k}\Phi_k\big{)} X^{t,x,*}_{k} +\big{(}C_{k,k}+D_{k,k}\Phi_k\big{)}X^{t,x,*}_kw_k, \\ [1mm]
X^{t,x,*}_{t} = x,~~k\in \mathbb{T}_t,
\end{array}\right. \\
&&\left\{\begin{array}{l}
X^{k,\Phi}_{\ell+1} =\big{(}A_{k,\ell}+B_{k,\ell}\Phi_\ell\big{)} X^{k,\Phi}_{\ell}+\big{(}C_{k,\ell}+D_{k,\ell}\Phi_\ell\big{)}X^{k,\Phi}_\ell w_\ell, \\ [1mm]
X^{k,\Phi}_{k} = X^{t,x,*}_k,~~\ell\in \mathbb{T}_k,
\end{array}\right.\\
&&\label{equilibrium-state-closed-3}\left\{\begin{array}{l}
X^{k,u_k,\Phi}_{\ell+1} =\big{(}A_{k,\ell}+B_{k,\ell}\Phi_\ell\big{)} X^{k,u_k,\Phi}_{\ell}+\big{(}C_{k,\ell}+D_{k,\ell}\Phi_\ell\big{)}X^{k,u_k,\Phi}_\ell w_\ell,~~\ell\in \mathbb{T}_{k+1},\\[1mm]
X^{k,u_k,\Phi}_{k+1} =A_{k,k}X^{k,u_k,\Phi}_k+B_{k,k}u_k+\big{(}C_{k,k}X^{k,u_k,\Phi}_k+D_{k,k}u_k\big{)}w_k,\\[1mm]
X^{k,u_k,\Phi}_{k} = X^{t,x,*}_k
\end{array}\right.
\end{eqnarray}
with
%
\begin{eqnarray*}
&&(\Phi X^{k,\Phi})|_{\mathbb{T}_k}=\big{(}\Phi_k X_k^{k,\Phi},\cdots, \Phi_{N-1} X_{N-1}^{k,\Phi}\big{)},\\
&&(\Phi X^{k,u_k,\Phi})|_{\mathbb{T}_{k+1}}=\big{(}\Phi_{k+1} X_{k+1}^{k,u_t,\Phi},\cdots, \Phi_{N-1} X_{N-1}^{k,u_k,\Phi}\big{)}.
\end{eqnarray*}
%
%
\end{definition}

From Definition \ref{definition-closed-loop}, we know that $\Phi|_{\mathbb{T}_t}$ (the restriction of $\Phi$ on $\mathbb{T}_t$)
is a linear feedback equilibrium strategy. Hence, $\Phi$ is time-consistent.
The following theorem is concerned with the existence of linear
feedback equilibrium strategy, which is parallel to Theorem
\ref{Theorem-Equivalentce-open-loop}.

\begin{theorem}\label{Theorem-Equivalentce-closed-loop}
The following statements are equivalent.

(i) There exists a linear feedback equilibrium strategy of Problem (LQ).

(ii) There exists a $\Phi=\{\Phi_0,\cdots, \Phi_{N-1}\}$ with
$\Phi_t\in \mathbb{R}^{m\times n}, t\in \mathbb{T}$, such that for
any initial pair $(t,x)\in \mathbb{T}\times \mathbb{R}^n$ and $k\in
\mathbb{T}_t$, the following FBS$\Delta$E has a solution
$(X^{k,\Phi},Z^{k,\Phi})$
\begin{eqnarray}\label{system-adjoint-closed}
\left\{\begin{array}{l}
X^{k,\Phi}_{\ell+1} =\big{(}A_{k,\ell}+B_{k,\ell}\Phi_\ell\big{)} X^{k,\Phi}_{\ell} +\big{(}C_{k,\ell}+D_{k,\ell}\Phi_\ell\big{)}X^{k,\Phi}_\ell w_\ell, \\ [1mm]
Z_\ell^{k,\Phi}=\big{(}A_{k,\ell}+B_{k,\ell}\Phi_\ell\big{)}^T\mathbb{E}(Z_{\ell+1}^{k,\Phi}|\mathcal{F}_{\ell-1}) +\big{(}C_{k,\ell}+D_{k,\ell}\Phi_\ell\big{)}^T\mathbb{E}(Z_{\ell+1}^{k,\Phi}w_\ell|\mathcal{F}_{\ell-1})\\[1mm]
\hphantom{Z_\ell^{k,\Phi}=}+\big{(}\Phi_\ell^TR_{k,\ell}\Phi_\ell+Q_{k,\ell}\big{)}X_\ell^{k,\Phi}, \\ [1mm]
X^{k,\Phi}_k=X^{t,x,*}_k,~~Z_N^{k,\Phi}=G_kX^{k,\Phi}_N,~~\ell\in \mathbb{T}_k
\end{array}\right.
\end{eqnarray}
with the stationary condition
\begin{eqnarray}\label{stationary-condition-closed}
0=R_{k,k}\Phi_kX_k^{k,\Phi}+B_{k,k}^T\mathbb{E}(Z_{k+1}^{k,\Phi}|\mathcal{F}_{k-1})+D_{k,k}^T\mathbb{E}(Z_{k+1}^{k,\Phi}w_k|\mathcal{F}_{k-1}),~~k\in \mathbb{T}_t,
\end{eqnarray}
and the convexity condition
\begin{eqnarray}\label{convex-closed}
&&\inf_{\bar{u}_k\in L^2_\mathcal{F}(k; \mathbb{R}^m)}\Big{\{}\mathbb{E}\big{[}\bar{u}_k^TR_{k,k}\bar{u}_k\big{]}+\sum_{\ell=k}^{N-1}\mathbb{E}\Big{[}({Y}_\ell^{k,\bar{u}_k, \Phi})^T\big{(}Q_{k,\ell}+\Phi_\ell^T R_{k,\ell}\Phi_\ell\big{)} {Y}_\ell^{k,\bar{u}_k, \Phi}\Big{]}\nonumber \\[1mm]
&&\hphantom{\inf_{\bar{u}_k\in L^2_\mathcal{F}(k; \mathbb{R}^m)}}+\mathbb{E}\Big{[}({Y}_N^{k,\bar{u}_k, \Phi})^TG_{k} {Y}_N^{k,\bar{u}_k, \Phi}\Big{]} \Big{\}}\geq 0.
\end{eqnarray}
%
%
In (\ref{system-adjoint-closed}) (\ref{convex-closed}), $X^{t,x,*}, Y^{k,\bar{u}_k,\Phi}$ are given by (\ref{equilibrium-state-closed}) and
\begin{eqnarray}\label{system-y-closed-0}
\left\{\begin{array}{l}
{Y}^{k,\bar{u}_k,\Phi}_{\ell+1}=\big{(}A_{k,\ell}+B_{k,\ell}\Phi_\ell\big{)}{Y}^{k,\bar{u}_k,\Phi}_\ell+\big{(}C_{k,\ell}+D_{k,\ell}\Phi_\ell\big{)}{Y}^{k,\bar{u}_k,\Phi}_\ell w_\ell,~~~\ell\in \mathbb{T}_{k+1},\\[1mm]
Y^{k,\bar{u}_k,\Phi}_{k+1}=B_{k,k}\bar{u}_k+D_{k,k}\bar{u}_kw_k,\\[1mm]
{Y}^{k,\bar{u}_k,\Phi}_k=0.
\end{array}\right.
\end{eqnarray}
\end{theorem}

\emph{Proof}. See Appendix B. \hfill $\square$

Based on Theorem \ref{Theorem-Equivalentce-closed-loop}, the relationship between the existence of linear feedback equilibrium strategy and the solvability of a set of difference equations is established. It is stated in the following theorem.

\begin{theorem}\label{Theorem-Equivalentce-closed-loop-2}
The following statements are equivalent.

(i) Problem (LQ) admits a linear feedback equilibrium strategy.

(ii) The following set of equations
\begin{eqnarray}\label{GDRE-system-closed}
\left\{
\begin{array}{l}
\left\{\begin{array}{l}
\widetilde{P}_{t,k}=Q_{t,k}+\Phi_k^T R_{t,k}\Phi_k+(A_{t,k} +B_{t,k}\Phi_k)^T\widetilde{P}_{t,k+1}(A_{t,k}+B_{t,k}\Phi_k) \\ [1mm]
\hphantom{\widetilde{P}_{t,k}=}+(C_{t,k}+D_{t,k}\Phi_k)^T\widetilde{P}_{t,k+1}(C_{t,k}+D_{t,k}\Phi_k) \\ [1mm]
\widetilde{P}_{t,N}=G_t,\\[1mm]
k\in \mathbb{T}_t,
\end{array}\right. \\ [1mm]
{\widetilde{W}}_{t,t}{\widetilde{W}}_{t,t}^\dag{\widetilde{H}}_{t,t}-{\widetilde{H}}_{t,t}=0,\\[1mm]
\widetilde{W}_{t,t}\geq 0,\\[1mm]
t\in \mathbb{T}
\end{array}\right.
\end{eqnarray}
is solvable in the sense that the constrained conditions ${\widetilde{W}}_{t,t}{\widetilde{W}}_{t,t}^\dag{\widetilde{H}}_{t,t}-{\widetilde{H}}_{t,t}=0$ and $\widetilde{W}_{t,t}\geq 0, \forall t\in\mathbb{T}$, are satisfied. Here,
\begin{eqnarray}\label{GDRE-W-H-closed}
\left\{\begin{array}{l}
{\widetilde{W}}_{t,t}=R_{t,t}+B^T_{t,t}\widetilde{P}_{t,t+1}B_{t,t}+D^T_{t,t}\widetilde{P}_{t,t+1}D_{t,t},\\[1mm]
{\widetilde{H}}_{t,t}=B^T_{t,t}\widetilde{P}_{t,t+1}A_{t,t}+D^T_{t,t}\widetilde{P}_{t,t+1}C_{t,t},\\[1mm]
\Phi_t=-\widetilde{W}_{t,t}^\dagger \widetilde{H}_{t,t},\\[1mm]
%
t\in \mathbb{T}.
\end{array}\right.
\end{eqnarray}

Under any of above conditions, $\Phi$ given in
(\ref{GDRE-W-H-closed}) is a linear feedback equilibrium strategy.
\end{theorem}

\emph{Proof}. (i)$\Rightarrow$(ii). Let $\Phi$ be a linear feedback
equilibrium strategy. Then, for any initial pair $(t,x)$,
(\ref{system-adjoint-closed}) admits a solution.
Let $k=t$ in (\ref{system-adjoint-closed}) and (\ref{stationary-condition-closed}); noting $Z_N^{t,\Phi}=G_tX^{t,\Phi}_N$ and by results in \cite{Ni-Li-Zhang}, one can get
\begin{eqnarray}\label{decouple-closed}
Z^{t,\Phi}_{k}=\widetilde{P}_{t,k}X^{t,\Phi}_k, ~k\in \mathbb{T}_t,
\end{eqnarray}
where $\widetilde{P}_{t,k}, k\in {\mathbb{T}}_t$, are deterministic matrices and determined below.
Then, from (\ref{stationary-condition-closed}) and Lemma \ref{Lemma-matrix-equation}, we have
\begin{eqnarray*}\label{stationary-condition-ell-2-closed}
&&0=R_{t,t}\Phi_tX_t^{t,\Phi}+B_{t,t}^T\mathbb{E}(\widetilde{P}_{t,t+1}X_{t+1}^{t,\Phi}|\mathcal{F}_{t-1})+D_{t,t}^T\mathbb{E}(\widetilde{P}_{t,t+1}X_{t+1}^{t,\Phi}w_t|\mathcal{F}_{t-1})\nonumber\\[1mm]
&&\hphantom{0}=\Big{[}R_{t,t}\Phi_t+B^T_{t,t}\widetilde{P}_{t,t+1}(A_{t,t}+B_{t,t}\Phi_t)+D^T_{t,t}\widetilde{P}_{t,t+1}(C_{t,t}+D_{t,t}\Phi_t)\Big{]}X_{t}^{t,\Phi}.
\end{eqnarray*}
As $x=X_{t}^{t,\Phi}$ can be arbitrarily selected, we have %
\begin{eqnarray*}
0=R_{t,t}\Phi_t+B^T_{t,t}\widetilde{P}_{t,t+1}(A_{t,t}+B_{t,t}\Phi_t)+D^T_{t,t}\widetilde{P}_{t,t+1}(C_{t,t}+D_{t,t}\Phi_t).
\end{eqnarray*}
From Lemma \ref{Lemma-matrix-equation}, it follows that
\begin{eqnarray*}
\Phi_t=-\widetilde{W}_{t,t}^\dagger \widetilde{H}_{t,t},
\end{eqnarray*}
and
\begin{eqnarray*}
{\widetilde{W}}_{t,t}{\widetilde{W}}_{t,t}^\dag{\widetilde{H}}_{t,t}-{\widetilde{H}}_{t,t}=0.
\end{eqnarray*}
Note that
\begin{eqnarray*}\label{Theorem-Equivalentce-closed-loop-s-g-3}
&&Z_k^{t,\Phi}=\Big{\{}Q_{t,k}+\Phi_k^T R_{t,k}\Phi_k+(A_{t,k} +B_{t,k}\Phi_k)^T\widetilde{P}_{t,k+1}(A_{t,k}+B_{t,k}\Phi_k)\nonumber\\
&&\hphantom{Z_k^{t,\Phi}=}+(C_{t,k}+D_{t,k}\Phi_k)^T\widetilde{P}_{t,k+1}(C_{t,k}+D_{t,k}\Phi_k)\Big{\}}X_k^{t,\Phi}.
\end{eqnarray*}
For $k\in \mathbb{T}_t$, let
\begin{eqnarray*}
&&\widetilde{P}_{t,k}=Q_{t,k}+\Phi_k^T R_{t,k}\Phi_k+(A_{t,k} +B_{t,k}\Phi_k)^T\widetilde{P}_{t,k+1}(A_{t,k}+B_{t,k}\Phi_k)\nonumber\\
&&\hphantom{P_{t,k}=}+(C_{t,k}+D_{t,k}\Phi_k)^T\widetilde{P}_{t,k+1}(C_{t,k}+D_{t,k}\Phi_k).
\end{eqnarray*}
Then, (\ref{decouple-closed}) holds.
Furthermore, by the convexity condition (\ref{convex-closed}), we
have for any $\bar{u}_t\in L^2_{\mathcal{F}}(t;\mathbb{R}^m)$
\begin{eqnarray*}\label{convex-closed-2}
&&0\leq \mathbb{E}\big{[} \bar{u}_t^TR_{t,t}\bar{u}_t\big{]}+\mathbb{E}\big{[}(Y_N^{t,\bar{u}_t,\Phi})^T\widetilde{P}_{t,N} {Y}_N^{t,\bar{u}_t,\Phi}\big{]}+\sum_{\ell=t}^{N-1}\mathbb{E}\Big{[}(Y_k^{t,\bar{u}_t, \Phi})^T\big{(}Q_{t,k}+\Phi_k^TR_{t,k}\Phi_k\big{)} {Y}_k^{t,\bar{u}_t,\Phi}\Big{]}\nonumber\\
&&\hphantom{0}=\mathbb{E}\big{[} \bar{u}_t^T\big{(}R_{t,t}+B_{t,t}^T\widetilde{P}_{t,t+1}B_{t,t}+D_{t,t}^T\widetilde{P}_{t,t+1}D_{t,t}\big{)}\bar{u}_t\big{]}.
\end{eqnarray*}
Hence,
\begin{eqnarray*}
\widetilde{W}_{t,t}=R_{t,t}+B_{t,t}^T\widetilde{P}_{t,t+1}B_{t,t}+D_{t,t}^T\widetilde{P}_{t,t+1}D_{t,t}\geq 0.
\end{eqnarray*}
When $t$ ranges from $N-1$ to $0$, we have the solvability of (\ref{GDRE-system-closed}).

(ii)$\Rightarrow$(i). Note that $\Phi=\{\Phi_0,...,\Phi_{N-1}\}$ with $\Phi_t=-\widetilde{W}_{t,t}^\dagger \widetilde{H}_{t,t}, t\in \mathbb{T}$.
As (\ref{GDRE-system-closed}) is solvable, FBS$\Delta$E (\ref{system-adjoint-closed}) is solvable with property (\ref{decouple-closed}). Furthermore, by reversing some presentations in ``(i)$\Rightarrow$(ii)", the stationary condition (\ref{stationary-condition-closed}) and the convexity condition (\ref{convex-closed}) are both satisfied. Therefore, $\Phi$ is a linear feedback equilibrium strategy.
\hfill $\square$

Substituting $\Phi$ into (\ref{GDRE-system-closed}), one get
\begin{eqnarray}\label{GDRE-system-closed-2}
\left\{\begin{array}{l}
\left\{\begin{array}{l}
\widetilde{P}_{t,k}=Q_{t,k}+A_{t,k}^T\widetilde{P}_{t,k+1}A_{t,k}+C_{t,k}^T\widetilde{P}_{t,k+1}C_{t,k}-\widetilde{H}^T_{k,k}\widetilde{W}^\dagger_{k,k} \widetilde{H}_{t,k}\\[1mm]
\hphantom{\widetilde{P}_{t,k}=}-\widetilde{H}_{t,k}^T\widetilde{W}^\dagger_{k,k}\widetilde{H}_{k,k}+\widetilde{H}_{k,k}^T\widetilde{W}^\dagger_{k,k} \widetilde{W}_{t,k}\widetilde{W}^{\dagger}_{k,k}\widetilde{H}_{k,k}\\[1mm]
\widetilde{P}_{t,N}=G_t,\\[1mm]
k\in \mathbb{T}_t,
\end{array}\right. \\ [1mm]
{\widetilde{W}}_{t,t}{\widetilde{W}}_{t,t}^\dag{\widetilde{H}}_{t,t}-{\widetilde{H}}_{t,t}=0,\\[1mm]
\widetilde{W}_{t,t}\geq 0, \\ [1mm]
t\in \mathbb{T},
\end{array}\right.
\end{eqnarray}
where
\begin{eqnarray*}
\left\{\begin{array}{l}
{\widetilde{W}}_{t,k}=R_{t,k}+B^T_{t,k}\widetilde{P}_{t,k+1}B_{t,k}+D^T_{t,k}\widetilde{P}_{t,k+1}D_{t,k},\\[1mm]
{\widetilde{H}}_{t,k}=B^T_{t,k}\widetilde{P}_{t,k+1}A_{t,k}+D^T_{t,k}\widetilde{P}_{t,k+1}C_{t,k},\\[1mm]
%
k\in \mathbb{T}_t,~~t\in \mathbb{T}.
\end{array}\right.
\end{eqnarray*}
For  any $t\in \mathbb{T}, k\in \mathbb{T}_t$, $\widetilde{P}_{t,k}$ is symmetric.
Furthermore, if $Q_{t,k}\geq 0, R_{t,k}> 0 , G_t\geq 0$, $k\in \mathbb{T}_t, t\in \mathbb{T}$, then we can prove that $\widetilde{W}_{t,k}>0$ and $\widetilde{P}_{t,k}\geq 0, k\in \mathbb{T}_t, t\in \mathbb{T}$. Therefore, (\ref{GDRE-system-closed-2}) is solvable.

\begin{corollary}
If $Q_{t,k}\geq 0, R_{t,k}> 0 , G_t\geq 0$, $k\in \mathbb{T}_t, t\in \mathbb{T}$, then Problem (LQ) admits a unique linear feedback equilibrium strategy $\Phi$.
\end{corollary}
%

%
For the case where $A_{t,k}, B_{t,k}, C_{t,k}, D_{t,k}, Q_{t,k}, k\in \mathbb{T}_t, G_t$, are all independent of $t$, (\ref{GDRE-system-closed-2}) reads as
\begin{eqnarray}\label{GDRE-system-closed-3}
\left\{
\begin{array}{l}
\left\{
\begin{array}{l}
\widetilde{P}_{t,k}=Q_{k}+A_{k}^T\widetilde{P}_{t,k+1}A_{k}+C_{k}^T\widetilde{P}_{t,k+1}C_{k}-\widetilde{H}^T_{k,k}\widetilde{W}^\dagger_{k,k} \widetilde{H}_{t,k}\\[1mm]
\hphantom{\widetilde{P}_{t,k}=}-\widetilde{H}_{t,k}^T\widetilde{W}^\dagger_{k,k}\widetilde{H}_{k,k}+\widetilde{H}_{k,k}^T\widetilde{W}^\dagger_{k,k} \widetilde{W}_{t,k}\widetilde{W}^{\dagger}_{k,k}\widetilde{H}_{k,k}\\[1mm]
\widetilde{P}_{t,N}=G,\\[1mm]
k\in \mathbb{T}_t,
\end{array}
\right.\\[1mm]
{\widetilde{W}}_{t,t}{\widetilde{W}}_{t,t}^\dag{\widetilde{H}}_{t,t}-{\widetilde{H}}_{t,t}=0,\\[1mm]
\widetilde{W}_{t,t}\geq 0,\\[1mm]
t\in \mathbb{T},
\end{array}\right.
\end{eqnarray}
where
\begin{eqnarray*}
\left\{\begin{array}{l}
{\widetilde{W}}_{t,k}=R_{t,k}+B^T_{k}\widetilde{P}_{t,k+1}B_{k}+D^T_{k}\widetilde{P}_{t,k+1}D_{k},\\[1mm]
{\widetilde{H}}_{t,k}=B^T_{k}\widetilde{P}_{t,k+1}A_{k}+D^T_{k}\widetilde{P}_{t,k+1}C_{k},\\[1mm]
%
k\in \mathbb{T}_t,~~t\in \mathbb{T}.
\end{array}\right.
\end{eqnarray*}
Here, (\ref{GDRE-system-closed-3}) is different from
(\ref{GDRE-system2-s}) and (\ref{LRE-s}), which relates to the
open-loop equilibrium control of the corresponding situation.
Moreover, when all the system matrices in the dynamics and cost
functional are independent of the initial time,
(\ref{GDRE-system-closed-2}) will reduce to the standard GDRE
\cite{Ait-Chen-Zhou-2002}.
%

\section{Comparison}

For the time-inconsistent stochastic LQ problem, open-loop and
closed-loop time-consistent solutions are separately investigated in
the above two sections. Section \ref{Section--open-loop} is
concerned with the open-loop equilibrium control, while the linear
feedback equilibrium strategy is studied in Section
\ref{Section--closed-loop} that is a kind of closed-loop strategy.
As noted in Introduction, the objects of study in these two
formulations are quite  different. Open-loop time-consistent
solution is to find an open-loop control that is an equilibrium of a
leader-follower game with hierarchical structure. Open-loop control,
or simply control, or more exactly control action, in this paper is
referred to  as a function of time that is also adapted to a
filtration; such a meaning of open-loop control is adopted by many
scholars in their works (for example, \cite{Basar}
\cite{Sun-jingrui}). Generated by some {primitive random variables},
the filtration is viewed as the ``\emph{state of nature}",  and the
adaptedness to such a filtration is to emphasize that the
underlining open-loop control is allowed to be random. Note that the
key point of above comment is not to the attribute ``open-loop" but
to the subject ``control".

So far, time-consistent (linear feedback) strategy is the object of
closed-loop formulation. A strategy is a decision rule that a
controller uses to select his control action based on the available
information set. Mathematically, strategy is a mapping or operator
on some information set, which is a higher-rank notion other than
control. When substituting the available information into a
strategy, the open-loop value or open-loop realization of this
strategy is then obtained. To get more about this, let us look at
the strategy $\Phi$ in the definition of linear feedback equilibrium
strategy, which can be viewed as a binary mapping $\Phi(\cdot,
\cdot)$. When substituting the information $(k, X^{t,x,*}_k)$, we
have the open-loop value $\Phi(k, X^{t,x,*}_k)=\Phi_k X^{t,x,*}_k$
at time point $k$. Similarly, $(\Phi
X^{k,\Phi})|_{\mathbb{T}_k}=(\Phi(k, X^{k,\Phi}_k),...,\Phi(N-1,
X^{k,\Phi}_{N-1}))$ is an open-loop value of $\Phi$ on the time
period $\mathbb{T}_k$, which is a control action on $\mathbb{T}_k$.

Besides above conceptual difference between control and strategy, let us further compare the definitions of the two equilibria. In the definition of open-loop equilibrium control, state involved in $J(k,X^{t,x,*}_k;(u_k, u^{t,x,*}|_{\mathbb{T}_{k+1}}))$ of (\ref{open-loop-equilibrium}) satisfies the following equations
\begin{eqnarray}\label{comparision-1}
\left\{\begin{array}{l}
X^{k,u_k,u^{t,x,*}}_{\ell+1} =A_{k,\ell}X^{k,u_k,u^{t,x,*}}_\ell+B_{k,\ell}u^{t,x,*}_\ell+\big{(}C_{k,\ell}X^{k,u_k,u^{t,x,*}}_\ell+D_{k,\ell}u^{t,x,*}_\ell\big{)}w_\ell,~~\ell\in \mathbb{T}_{k+1},\\[1mm]
X^{k,u_k,u^{t,x,*}}_{k+1} =A_{k,k}X^{k,u_k,u^{t,x,*}}_k+B_{k,k}u_k+\big{(}C_{k,k}X^{k,u_k,u^{t,x,*}}_k+D_{k,k}u_k\big{)}w_k,\\[1mm]
X^{k,u_k,u^{t,x,*}}_{k} = X^{t,x,*}_k.
\end{array}\right.
\end{eqnarray}
In (\ref{comparision-1}), $u^{t,x,*}|_{\mathbb{T}_{k+1}}$ is not influenced by $u_k$, because it is given prior.
On the contrary, in the definition of linear feedback equilibrium
strategy, the control $(\Phi X^{k,u_k,\Phi})|_{\mathbb{T}_{k+1}}$ in
(\ref{equilibrium-state-closed-3}) is influenced by ${u}_k$ via the
term $X^{k,u_k,\Phi}$. This makes an essential difference between
the open-loop equilibrium control and the linear feedback
equilibrium strategy. In addition, an open-loop equilibrium control
when being mentioned is corresponding to a fixed initial pair, while
the linear feedback equilibrium strategy is required to define for
all the initial pairs.

Moreover, the existence of open-loop equilibrium control and linear
feedback equilibrium strategy is differently characterized. The fact
that there exists an  open-loop equilibrium control for any given
initial pair, is fully characterized via the solvability of a set of
nonsymmetric GDREs and the solvability of a set of LDEs. In
contrast, the existence of linear feedback equilibrium strategy is
shown to be equivalent to the solvability of a single set of GDREs,
which are symmetric. Note that just $R_{k,k}, k\in \mathbb{T}$, are
involved in (\ref{GDRE-system2}) and (\ref{LRE}). Hence, if some
$R_{t,k}, k\in \mathbb{T}_{t+1}, t\in \mathbb{T}$, are modified, the
existence of open-loop equilibrium control will not change! The
reason of this lies in the definition of open-loop equilibrium
control. Specifically, given initial pair $(t,x)$ and by subtracting
$
\sum_{\ell=k+1}^{N-1}\mathbb{E}\big{[}(u^{t,x,*}_\ell)^TR_{k,\ell}u^{t,x,*}_\ell \big{]}
$
from both sides of  (\ref{open-loop-equilibrium}), (\ref{open-loop-equilibrium}) can be equivalently rewritten as
\begin{eqnarray}\label{comparision-2}
&&\hspace{-3em}\mathbb{E}\Big{[}\sum_{\ell=k}^{N-1}(X^{k,u^{t,x,*}}_\ell)^TQ_{k,\ell}X^{k,u^{t,x,*}}_\ell \Big{]}+\mathbb{E}\big{[}(u^{t,x,*}_k)^TR_{k,k}u^{t,x,*}_k \big{]}+\mathbb{E}\big{[}(X^{k,u^{t,x,*}}_N)^TG_{k}X^{k,u^{t,x,*}}_N \big{]}\nonumber \\
&&\hspace{-3em}\leq \mathbb{E}\Big{[}\sum_{\ell=k}^{N-1}(X^{k,u_k,u^{t,x,*}}_\ell)^TQ_{k,\ell}X^{k,u_k,u^{t,x,*}}_\ell \Big{]}+\mathbb{E}\big{[}u_k^TR_{k,k}u_k \big{]}+\mathbb{E}\big{[}(X^{k,u_k,u^{t,x,*}}_N)^TG_{k}X^{k,u_k,u^{t,x,*}}_N \big{]},
\end{eqnarray}
where $X^{k,u^{t,x,*}}$ and $X^{k,u_k,u^{t,x,*}}$ are given, respectively, by
\begin{eqnarray}\label{comparision-3}
\left\{\begin{array}{l}
X^{k,u^{t,x,*}}_{\ell+1} =A_{k,\ell}X^{k,u^{t,x,*}}_\ell+B_{k,\ell}u^{t,x,*}_\ell+\big{(}C_{k,\ell}X^{k,u^{t,x,*}}_\ell+D_{k,\ell}u^{t,x,*}_\ell\big{)}w_\ell,\\[1mm]
X^{k,u^{t,x,*}}_{k} = X^{t,x,*}_k,~~\ell\in \mathbb{T}_{k+1},
\end{array}\right.
\end{eqnarray}
and (\ref{comparision-1}). Note that $R_{k,\ell}, \ell\in
\mathbb{T}_{k+1}, k\in \mathbb{T}$, do not appear in
(\ref{comparision-2}). This could explain why only $R_{k,k}, k\in
\mathbb{T}$, are involved in (\ref{GDRE-system2}) and (\ref{LRE}).
Otherwise, in the definition of linear feedback equilibrium
strategy, $(\Phi X^{k,\Phi})|_{\mathbb{T}_{k+1}}$ and $(\Phi
X^{k,u_k,\Phi})|_{\mathbb{T}_{k+1}}$ of (\ref{defi-closed-loop}) are
different as $X^{k,\Phi}$ differs from $X^{k,u_k,\Phi}$. Therefore,
terms in (\ref{defi-closed-loop}) associated with $R_{k,\ell},
\ell\in \mathbb{T}_{k+1}$, cannot be removed, which implies the
dependence of (\ref{GDRE-system-closed-2}) on $R_{k,\ell}, \ell\in
\mathbb{T}_{k+1}, k\in \mathbb{T}$.

Furthermore, let us pay attention to the computations of these two
equilibria. As noted above, when we mention an open-loop equilibrium
control, the initial pair which induces that open-loop equilibrium
control should be mentioned simultaneously. Differently, the linear
feedback equilibrium strategy is independent of all the initial
pairs. These facts will help us to differentiate the backward
procedures (mentioned below) of computing an open-loop equilibrium
control and the linear feedback equilibrium strategy. Theorem
\ref{Theorem-Equivalentce-open-loop} and Theorem
\ref{Theorem-Equivalentce-closed-loop} are generally the
maximum-principle-type equivalent characterizations, and the
backward S$\Delta$Es  are involved. To obtain an open-loop
equilibrium control $u^{t,x,*}$, we should decouple the
FBS$\Delta$Es (for every $k$ we have a (\ref{system-adjoint})) along
the equilibrium state $X^{t,x,*}$. Roughly speaking, from the
stationary condition,  we have $u^{t,x,*}_{N-1}$ by decoupling
FBS$\Delta$E (\ref{system-adjoint}) (with $k=N-1$).
Generally, based on all the expressions of $u^{t,x,*}_\ell, \ell\in \mathbb{T}_{k+1}$, we should redecouple the FBS$\Delta$E (\ref{system-adjoint}) on the whole time period $\mathbb{T}_k$ to obtain the linear relation between $Z^{k,*}$ and $X^{k,*}$, and then by the stationary condition we can get $u^{t,x,*}_k$. %
Note that via its forward initial state $X^{t,x,*}_k$, the FBS$\Delta$Es are coupled with the equilibrium state $X^{t,x,*}$, and the equilibrium control and the equilibrium state are obtained simultaneously.

Concerned with the linear feedback equilibrium strategy, if exists, it can be calculated via solving the GDREs (\ref{GDRE-system-closed-2}) backwardly. For every generic  time point $t$, we need to solve a single GDRE (associated with $t$) over the whole time period $\mathbb{T}_t$, and use the solution $\widetilde{P}_{t,t+1}$ at time point $t+1$ to construct $\Phi_t=\widetilde{W}_{t,t}^\dagger \widetilde{H}_{t,t}$. This backward procedure is also due to a decoupling procedure of FBS$\Delta$Es (\ref{system-adjoint-closed}).
It should be mention that only under the condition of the existence of linear feedback equilibrium strategy, this strategy can be computed backwardly as above.
Furthermore, necessary and sufficient conditions on the existence of
the equilibria are the main concerns of this paper, from which the
computations of the equilibria are much direct.

To end this section, we give the following final comment. Under the
condition that there exists an open-loop equilibrium control for
any given initial pair, all the open-loop equilibrium controls
happen to be of feedback form; while the feedback gain
$\{-W^\dagger_{k,k} H_{k,k}, k\in \mathbb{T}\}$ in
(\ref{Theorem-Equivalentce-open-loop-control}) is not a linear
feedback equilibrium strategy of Problem (LQ) indeed!

\section{Examples}

In this section, we shall present three examples to illustrate the theory derived above.

\begin{example}\label{Example---1}
Let
\begin{eqnarray*}
&&A_0=\left[
\begin{array}{cc}
1.12 & 0.21\\
-0.13& 0.98
\end{array}
\right],~~
A_1=\left[
\begin{array}{cc}
2.12 & -0.35\\
-0.21& 3.43
\end{array}
\right],~~A_2=\left[
\begin{array}{cc}
5.46 & 1.21\\
-0.98& 4.21
\end{array}
\right],\\
&&B_0=\left[
\begin{array}{cc}
1.45 & -0.23\\
-0.2& 4
\end{array}
\right],~~B_1=\left[
\begin{array}{cc}
1.5 & 0.3\\
-0.2& 3
\end{array}
\right],~~B_2=\left[
\begin{array}{cc}
-4.36 & 0.82\\
1.21& 4.21
\end{array}
\right],\\[1mm]
&&C_0=\left[
\begin{array}{cc}
1 & 0.32\\
0.25&  3
\end{array}
\right],~~C_1=\left[
\begin{array}{cc}
1.65 & -0.13\\
-0.42& 6
\end{array}
\right],~~C_2=\left[
\begin{array}{cc}
-3 & 1.53\\
-0.62& 4.78
\end{array}
\right],\\[1mm]
&&D_0=\left[
\begin{array}{cc}
5 & 1\\
-0.85&  8
\end{array}
\right],~~
D_1=\left[
\begin{array}{cc}
4 & 0.53\\
-0.42& 5
\end{array}
\right],~~D_2=\left[
\begin{array}{cc}
9.21 & -2.03\\
-1.52& 6.98
\end{array}
\right],\\[1mm]
&&Q_0=\left[
\begin{array}{cc}
-2 & 0.8\\
0.8&  1.6
\end{array}
\right],~~Q_1=\left[
\begin{array}{cc}
4 & 0\\
0& 0
\end{array}
\right],~~Q_2=\left[
\begin{array}{cc}
1.56 & -0.23\\
-0.23& 2.54
\end{array}
\right],~~R_{0,0}=\left[
\begin{array}{cc}
-0.5 & 0\\
0&  1
\end{array}
\right],\\
&&R_{0,1}=\left[
\begin{array}{cc}
-5 & 0\\
0&  -4
\end{array}
\right],~~R_{0,2}=\left[
\begin{array}{cc}
-9 & 0\\
0&  10
\end{array}
\right],~~R_{1,1}=\left[
\begin{array}{cc}
4 & -0.3\\
-0.3& 7
\end{array}
\right],\\[1mm]
&&R_{1,2}=\left[
\begin{array}{cc}
2.24 & -5.67\\
-5.67&  -1.27
\end{array}
\right],~~R_{2,2}=\left[
\begin{array}{cc}
6.29 & -1.67\\
-1.67&  8.38
\end{array}
\right],~~G=\left[
\begin{array}{cc}
1 & 0\\
0&  2
\end{array}
\right].
\end{eqnarray*}
Here, $Q_1\geq 0, Q_2, G,  R_{1,1}, R_{2,2}>0$; $Q_0, R_{0,0}$ and $R_{1,2}$ are indefinite, and $R_{0,1}$ is negative definite.\end{example}

By (\ref{GDRE-system2-s}) and (\ref{LRE-s}), we get
\begin{eqnarray*}
&&\hspace{-1.5em}P_2=\left[
\begin{array}{cc}
   16.6571 &   5.8520\\
    5.8520 &  11.5436
\end{array}
\right],  ~~
P_1=\left[
\begin{array}{cc}
    6.9700 &  -1.3882\\
   -1.3882 &   9.1396
\end{array}
\right],~~P_0=\left[
\begin{array}{cc}
    7.8991 &   4.2276\\
    4.2276 &   4.6336
\end{array}
\right],\\
&&\hspace{-1.5em}S_2=\left[\begin{array}{cc}
   43.0612 & -12.3922\\
  -12.3922 &  87.4900
\end{array}
\right],~~S_1=\left[
\begin{array}{cc}
   11.6579 &  -7.4371\\
   -7.4371 &  95.6692
\end{array}
\right],~~
S_0=\left[
\begin{array}{cc}
   34.3248 &  35.9699\\
   35.9699 & 938.8710
\end{array}
\right]
\end{eqnarray*}
with
\begin{eqnarray}
&&\label{Example-4.1-2} W_2=\left[
\begin{array}{cc}
 117.6727 & -34.9725\\
  -34.9725&  146.0623
\end{array}
\right]>0,~~W_1=\left[
\begin{array}{cc}
  287.3160 & 153.0623\\
  153.0623 & 447.2115
\end{array}
\right]>0,\\[1mm]
&&\label{Example-4.1-0}W_0=\left[
\begin{array}{cc}
    1138.3 &  -410.9\\
   -410.9  &  915.8
\end{array}
\right]>0,\\
&&R_{2,2}+B^T_{2}S_{3}B_{2}+D^T_{2}S_{3}D_{2}=\left[
\begin{array}{cc}
  117.6727 & -34.9725\\
  -34.9725 & 146.0623
\end{array}
\right]>0,\\[1mm]
&&R_{1,1}+B^T_{1}S_{2}B_{1}+D^T_{1}S_{2}D_{1}=\left[
\begin{array}{cc}
    857.9 &  -426.0\\
   -426.0 &   2909.6
\end{array}
\right]>0,\\[1mm]
&&\label{Example-4.1-4}R_{0,0}+B^T_{0}S_{1}B_{0}+D^T_{0}S_{1}D_{0}=\left[
\begin{array}{cc}
    17940  &-54740\\
   -54740  &  331470
\end{array}
\right]>0.
\end{eqnarray}
From (\ref{Example-4.1-2})-(\ref{Example-4.1-4}), the corresponding (\ref{GDRE-system2-s}) and (\ref{LRE-s}) are solvable, and thus, the open-loop equilibrium pair exists. Furthermore, an open-loop equilibrium control for the initial pair $(0,x)$ is given by
\begin{eqnarray*}
u^{0,x,*}_k=-W_k^\dagger H_kX^{0,x,*}_k, ~~k=0,1,2
\end{eqnarray*}
with
\begin{eqnarray*}
&&-W_0^\dagger H_0=\left[
\begin{array}{cc}
   -0.2183  &  0.0031\\
    0.0023  & -0.3286
\end{array}
\right],~~
-W_1^\dagger H_1=\left[
\begin{array}{cc}
   -0.5138  &  0.1973\\
    0.0026  & -1.1339
\end{array}
\right],\\
&&
-W_2^\dagger H_2=\left[
\begin{array}{cc}
    0.4889 &  -0.2601\\
    0.1605 &  -0.7474
\end{array}
\right]
\end{eqnarray*}
and
\begin{eqnarray*}
\left\{
\begin{array}{l}
X^{0,x,*}_{k+1} = (A_{k}X^{0,x,*}_{k}+B_{k}u^{0,x,*}_k)\\
\hphantom{X^{0,x,*}_{k+1} = }+(C_{k}X^{0,x,*}_{k} +D_{k}u^{0,x,*}_k)w_k,\\
X^{0,x,*}_{0} = x,~~ k=0,1,2.
\end{array}
\right.
\end{eqnarray*}

On the other hand, from (\ref{GDRE-system-closed-3}) we can obtain the solution.
However, we have
\begin{eqnarray*}
\widetilde{W}_{1,1}=\left[
\begin{array}{cc}
  239.0218 & 247.7565\\
  247.7565 & 224.5117
   \end{array}
   \right],
\end{eqnarray*}
whose eigenvalues are
\begin{eqnarray*}
\lambda_1=-16.096,~~~\lambda_2= 479.6294.
\end{eqnarray*}
Clearly, $\widetilde{W}_{1,1}$ is indefinite, and thus, the corresponding (\ref{GDRE-system-closed-3}) is not solvable.
This means that the linear feedback equilibrium strategy does not exist.

\begin{example}\label{Example---2}
The system matrices and the weight matrices are the same as those of
Example \ref{Example---1} except for $R_{0,1}, R_{0,2}, R_{1,2}$,
which are now
\begin{eqnarray*}
R_{0,1}=\left[
\begin{array}{cc}
-1 & 0\\
0&  -0.6
\end{array}
\right],~~R_{0,2}=\left[
\begin{array}{cc}
9.45 & 1.32\\
1.32&  10.78
\end{array}
\right],~~R_{1,2}=\left[
\begin{array}{cc}
5.24 & -1.67\\
-1.67&  7.27
\end{array}
\right].
\end{eqnarray*}
Note that $R_{0,1}$ is negative definite.
\end{example}

In this case, the corresponding  (\ref{GDRE-system2-s}) and
(\ref{LRE-s}) are the same as those of Example \ref{Example---1},
because $R_{0,1}, R_{0,2}, R_{1,2}$ do not enter the GDRE and LDE.
Hence, the open-loop equilibrium pair exists for any initial pair.
Let us check the existence of the linear feedback equilibrium
strategy. By (\ref{GDRE-system-closed-3}) we have
\begin{eqnarray*}
&&\widetilde{P}_{2,2}=\left[\begin{array}{cc}
   16.6571 &   5.8520\\
    5.8520 &  11.5436
   \end{array}
   \right],~~\widetilde{P}_{1,2}=\left[\begin{array}{cc}
   16.3775 &   6.1187\\
    6.1187 &  10.8526
   \end{array}
   \right],~~\widetilde{P}_{1,1}=\left[\begin{array}{cc}
   37.3769 & -10.7301\\
  -10.7301 &  12.7823
   \end{array}
   \right],\\[1mm]
&&\widetilde{P}_{0,2}=\left[\begin{array}{cc}
   17.9435  &  3.9449\\
    3.9449  & 14.2605
   \end{array}
   \right],~~\widetilde{P}_{0,1}=\left[\begin{array}{cc}
   39.7057 & -11.0096\\
  -11.0096 &   3.2054
   \end{array}
   \right],~~\widetilde{P}_{0,0}=\left[\begin{array}{cc}
    6.1615 &   4.3853\\
    4.3853 &   3.2889
   \end{array}
   \right]
\end{eqnarray*}
with
\begin{eqnarray}
&&\label{Example-4.2-2}\widetilde{W}_{2,2}=\left[\begin{array}{cc}
  117.6727 & -34.9725\\
  -34.9725 & 146.0623
   \end{array}
   \right]>0,~~\widetilde{W}_{1,1}=\left[\begin{array}{cc}
  281.0078 & 160.6675\\
  160.6675 & 425.5062
   \end{array}
   \right]>0,\\[1mm]
&&\label{Example-4.2-0}\widetilde{W}_{0,0}=\left[\begin{array}{cc}
    1178.0 &  -334.5\\
   -334.5  &  143.3
   \end{array}
   \right]>0.
\end{eqnarray}
It follows from (\ref{Example-4.2-2})-(\ref{Example-4.2-0}) that (\ref{GDRE-system-closed-3}) is solvable. Furthermore, a linear feedback equilibrium strategy $(\Phi_0, \Phi_1, \Phi_2)$ is given by
\begin{eqnarray*}
&&\Phi_0=\left[\begin{array}{cc}
   -0.0368  &  0.0884\\
    0.6555  & -0.0192
   \end{array}
   \right],~~\Phi_1=\left[\begin{array}{cc}
   -0.5094  &  0.1935\\
   -0.0021  & -1.1301
   \end{array}
   \right],~~\Phi_2=\left[\begin{array}{cc}
    0.4889 &  -0.2601\\
    0.1605 &  -0.7474
   \end{array}
   \right].
\end{eqnarray*}

\begin{example}
Let
\begin{eqnarray*}
&&A_{0,0}=\left[
\begin{array}{cc}
2.3 & 0.41\\
-0.3& 1.9
\end{array}
\right],~~
A_{0,1}=\left[
\begin{array}{cc}
4.12 & -0.35\\
0.31& 3.03
\end{array}
\right],~~B_{0,0}=\left[
\begin{array}{cc}
2.45 & -0.3\\
0.2& 4
\end{array}
\right],~~\\[1mm]
&&B_{0,1}=\left[
\begin{array}{cc}
2.5 & 0.6\\
-0.2& 3
\end{array}
\right],~~C_{0,0}=\left[
\begin{array}{cc}
2.2 & 0.32\\
0.5&  3
\end{array}
\right],~~C_{0,1}=\left[
\begin{array}{cc}
3.65 & -0.3\\
-0.42& 5.6
\end{array}
\right],\\[1mm]
&&D_{0,0}=\left[
\begin{array}{cc}
5.6 & 1\\
0.73&  7.8
\end{array}
\right],~~D_{0,1}=\left[
\begin{array}{cc}
5 & 0.73\\
-0.47& 5.2
\end{array}
\right],~~A_{1,1}=\left[
\begin{array}{cc}
6 & 1.63\\
-1.37& 7
\end{array}
\right],\\
&&B_{1,1}=\left[
\begin{array}{cc}
4 & 0.93\\
1.07& 3
\end{array}
\right],~~C_{1,1}=\left[
\begin{array}{cc}
8 & 2.03\\
-1.23& 10
\end{array}
\right],~~
D_{1,1}=\left[
\begin{array}{cc}
5 & -0.93\\
1.016 & 4.65
\end{array}
\right],\\[1mm]
&&Q_{0,0}=\left[
\begin{array}{cc}
2 & 0.8\\
0.8&  1.6
\end{array}
\right],~~Q_{0,1}=\left[
\begin{array}{cc}
4 & 0\\
0& 0
\end{array}
\right],~~R_{0,0}=\left[
\begin{array}{cc}
-0.5 & 0\\
0&  1
\end{array}
\right],~~R_{0,1}=\left[
\begin{array}{cc}
-5 & 0\\
0&  -4
\end{array}
\right],\\[1mm]
&&Q_{1,1}=\left[
\begin{array}{cc}
2 & 0.1\\
0.1& 5
\end{array}
\right],~~R_{1,1}=\left[
\begin{array}{cc}
4 & -0.3\\
-0.3& 7
\end{array}
\right],~~G_0=\left[
\begin{array}{cc}
1 & 0\\
0&  2
\end{array}
\right],~~G_1=\left[
\begin{array}{cc}
2 & -0.3\\
-0.3&  3
\end{array}
\right].
\end{eqnarray*}
Note that $R_{0,0}$ is indefinite and $R_{0,1}$ is negative definite.
\end{example}

By (\ref{GDRE-system-closed-2}), we have
\begin{eqnarray*}
\widetilde{P}_{1,1}=\left[\begin{array}{cc}
   18.8304 & -11.9513\\
  -11.9513 &  46.5418
   \end{array}
   \right],~~\widetilde{P}_{0,1}=\left[\begin{array}{cc}
   40.6027 & -28.7266\\
  -28.7266 &  50.9647
   \end{array}
   \right],~~\widetilde{P}_{0,0}=\left[\begin{array}{cc}
   99.6787 &  14.1112\\
   14.1112 &   8.3265
   \end{array}
   \right]
\end{eqnarray*}
with
\begin{eqnarray}
\label{Example-4.3-1}\widetilde{W}_{1,1}=\left[\begin{array}{cc}
   86.9155 &  11.0531\\
   11.0531 & 103.2478
   \end{array}
   \right]>0,~~\widetilde{W}_{0,0}=\left[\begin{array}{cc}
    1282.7  & -1027.0\\
   -1027.0  &  3582.2
   \end{array}
   \right]>0.
\end{eqnarray}
It follows from (\ref{Example-4.3-1}) that the corresponding (\ref{GDRE-system-closed-2}) is solvable. Hence, the linear feedback equilibrium strategy does exist. Furthermore, a linear feedback equilibrium strategy $(\Phi_0, \Phi_1)$ is given by
\begin{eqnarray*}
\Phi_0=\left[\begin{array}{cc}
   -0.4665  & -0.0206\\
    0.0269  & -0.3965
   \end{array}
   \right],~~\Phi_1=\left[\begin{array}{cc}
   -1.4499 &  -0.4726\\
    0.6369 &  -1.8700
   \end{array}
   \right].
\end{eqnarray*}

\section{Conclusion}

In this paper, we investigated the open-loop equilibrium control and
the linear feedback equilibrium strategy of the time-inconsistent
indefinite stochastic LQ optimal control. Necessary and sufficient
conditions are presented for these two cases, respectively.
Furthermore, the GDREs and LDEs are introduced to characterize the
linear feedback form of the open-loop equilibrium control and the
linear feedback equilibrium strategy. For future researches, we
would like to study the time-inconsistent problems for jump
parameter systems \cite{Caines}, Boolean networks \cite{Guo},
multi-agent systems \cite{Qi} \emph{etc.}, and extend the
methodology developed in this paper to other types of
time-inconsistency.





\section*{Appendix}

\subsection*{A. Proof of Theorem \ref{Theorem-Equivalentce-open-loop}}\label{appendix-A}         

Denote $J(k,X^{t,x,*}_k;(u_k, u^{t,x,*}|_{\mathbb{T}_{k+1}}))$ by $\bar{J}(k, X^{t,x,*}_k; u_k)$.  
Then, by (\ref{open-loop-equilibrium}), we have
\begin{eqnarray}\label{open-loop-equilibrium-4}
\bar{J}(k,X^{t,x,*}_k; u_k^{t,x,*})\leq \bar{J}(k,X^{t,x,*}_k; u_k).
\end{eqnarray}
%
%
%
This means that $u^{t,x,*}_k$ is an optimal control of the following nonstandard optimal control problem (denoted as Problem (LQ)$_k$):
\begin{eqnarray}
\left\{\begin{array}{l}
\mbox{Minimize }\bar{J}(k,X^{t,x,*}_k; u_k)\mbox{ over }L^2_\mathcal{F}(k; \mathbb{R}^m),\nonumber\\
\mbox{subject to}\nonumber\\
\label{system-u-u*}
\left\{\begin{array}{l}
X^k_{\ell+1} = A_{k,\ell}X^k_\ell+B_{k,\ell}u^{t,x,*}_\ell +\big{(}C_{k,\ell}X^k_\ell+D_{k,\ell}u^{t,x,*}_\ell\big{)}w_\ell, \\ [0.5mm]
X^k_{k+1} = (A_{k,k}X^k_k+B_{k,k}u_k)+(C_{k,k}X^k_k +D_{k,k}u_k)w_k, \\[0.5mm]
X^k_{k} = X_k^{t,x,*},~~ \ell\in\mathbb{T}_{k+1}.
\end{array}\right.
\end{array}\right.
\end{eqnarray}
Here, we call Problem (LQ)$_k$ a non-standard optimal control problem as $u^{t,x,*}|\mathbb{T}_{k+1}$ in the dynamics of $X^k$ is fixed, and we just select $u_k$ to minimize $\bar{J}(k,X^{t,x,*}_k; u_k)$.

To proceed, we introduce an inner product on $L^2_{\mathcal{F}}(\mathbb{T}_k; \mathbb{R}^p)$ with $p=n, m$, and $k\in \widetilde{\mathbb{T}}_t$:
\begin{eqnarray*}
\langle y, z\rangle_{\mathbb{T}_{k}}= \sum_{\ell=k}^{N-1}\mathbb{E}(y_\ell^Tz_\ell),\mbox{ for }y, z\in L^2_{\mathcal{F}}(\mathbb{T}_k; \mathbb{R}^p),
\end{eqnarray*}
and use the convention
\begin{eqnarray}\label{notation-1}
\left\{\begin{array}{l}
(Q_kx)(\cdot)=Q_{k,\cdot} x_\cdot,~~\forall x\in L^2_{\mathcal{F}}(\mathbb{T}_k; \mathbb{R}^n),\\[0.5mm]
%
(R_{k+}u)(\cdot)=R_{k,\cdot}u_\cdot,~~\forall u\in L^2_{\mathcal{F}}(\mathbb{T}_{k+1}; \mathbb{R}^m).
%
\end{array}\right.
\end{eqnarray}
For
$ L^2_\mathcal{F}(k; \mathbb{R}^p)$ with $p=n ,m$, and $k\in \widetilde{\mathbb{T}}_t$, the inner product
is defined as
\begin{eqnarray*}
\langle y, z\rangle_{k}=\mathbb{E}(y^Tz), \mbox{ for } y,z\in L^2_\mathcal{F}(k; \mathbb{R}^p).
\end{eqnarray*}
%
%
Then, the cost functional $\bar{J}(k, X^{t,x,*}_k; u_k)$ can be rewritten as
\begin{eqnarray}\label{cost-functional-2}
&&\bar{J}(k, X^{t,x,*}_k; u_k)=\langle Q_kX^k, X^k\rangle_{\mathbb{T}_k}+\langle R_{k,k}u_k, u_k\rangle_k+\langle R_{k+} u^{t,x,*}|_{\mathbb{T}_{k+1}}, u^{t,x,*}|_{\mathbb{T}_{k+1}}\rangle_{\mathbb{T}_{k+1}}\nonumber \\[1mm]
&&\hphantom{\bar{J}(k, X^{t,x,*}_k; u_k)=}+\langle G_kX^k_{N}, X^k_{N}\rangle_{N}.
\end{eqnarray}

We now calculate the first order and second order directional derivatives of $\bar{J}(k, X^{t,x,*}_k; u_k)$ at $u^{t,x,*}_k$. Corresponding to controls $u^{t,x,*}|_{\mathbb{T}_k}$ and $(u^{t,x,*}_k+\lambda\bar{u}_k, u^{t,x,*}|_{\mathbb{T}_{k+1}})$, the solutions of (\ref{system-k}) with the initial state $X^{t,x,*}_k$ are, respectively, denoted by $X^{k,*}$ and ${X}^{k,\lambda}$. Then, we have
\begin{eqnarray*}
\left\{\begin{array}{l}
\frac{{X}^{k,\lambda}_{\ell+1}-X^{k,*}_{\ell+1}}{\lambda}=A_{k,\ell}\frac{{X}_\ell^{k,\lambda}-X^{k,*}_\ell}{\lambda} +C_{k,\ell}\frac{{X}_\ell^{k,\lambda}-X^{k,*}_\ell}{\lambda}w_\ell, \\[1mm]
\frac{{X}^{k,\lambda}_{k+1}-X^{k,*}_{k+1}}{\lambda}=\Big{(}A_{k,k}\frac{{X}_k^{k,\lambda}-X^{k,*}_k}{\lambda}+B_{k,k}\bar{u}_k\Big{)}+\Big{(}C_{k,k}\frac{{X}_k^{k,\lambda}-X^{k,*}_k}{\lambda}+D_{k,k}\bar{u}_k\Big{)}w_k,\\
\frac{{X}_k^{k,\lambda}-X^{k,*}_k}{\lambda}=0,~~~\ell\in \mathbb{T}_{k+1},
\end{array}\right.
\end{eqnarray*}
which can be rewritten as
\begin{eqnarray}\label{system-y}
\left\{\begin{array}{l}
{Y}^k_{\ell+1}=A_{k,\ell}{Y}^k_\ell+C_{k,\ell}{Y}^k_\ell w_\ell,~~~\ell\in \mathbb{T}_{k+1},\\
Y^k_{k+1}=B_{k,k}\bar{u}_k+D_{k,k}\bar{u}_kw_k,\\
{Y}^k_k=0
\end{array}\right.
\end{eqnarray}
with $\frac{{X}_\ell^{k,\lambda}-X^{k,*}_\ell}{\lambda}={Y}^k_\ell$.
For any $\ell\in \mathbb{T}_k$, we have $X^{k,\lambda}_\ell=X^{k,*}_\ell+\lambda Y^k_\ell$.
To proceed, some calculations show
\begin{eqnarray*}\label{R}
&&\lim_{\lambda \downarrow 0}\frac{\langle R_{k,k}(u_k^{t,x,*}+\lambda \bar{u}_k),  u^{t,x,*}_k+\lambda\bar{u}_k\rangle_k-\langle R_{k,k}u_k^{t,x,*},  u_k^{t,x,*}\rangle_k}{\lambda} \\
&&=2\langle R_{k,k} u_k^{t,x,*}, \bar{u}_k\rangle_k+\lim_{\lambda\downarrow 0} \lambda\langle R_{k,k}\bar{u}_k, \bar{u}_k\rangle_t \nonumber \\
&&= 2\langle R_{k,k} u_k^{t,x,*}, \bar{u}_k\rangle_k,
\end{eqnarray*}
and
\begin{eqnarray*}
&&\label{Q}\lim_{\lambda \downarrow 0}\frac{\langle Q_{k}X^{k,\lambda}, X^{k,\lambda}\rangle_{\mathbb{T}_k}-\langle Q_{k}X^{k,*}, X^{k,*}\rangle_{\mathbb{T}_k}}{\lambda}= 2\langle Q_{k} X^{k,*}, Y^k\rangle_{\mathbb{T}_k},\\[1mm]
&&\label{G}\lim_{\lambda \downarrow 0}\frac{\langle G_{k}X_N^{k,\lambda}, X_N^{k,\lambda}\rangle_{N}-\langle G_{k}X_N^{k,*}, X_N^{k,*}\rangle_{N}}{\lambda}= 2\langle G_{k} X_N^{k,*}, Y_N^k\rangle_{N}.
\end{eqnarray*}
Hence, we have the first order directional derivative of $\bar{J}(k, X_k^{t,x,*};u_k)$ at $u_k^{t,x,*}$ with the direction $\bar{u}_k$
\begin{eqnarray}\label{first-order}
&&d\bar{J}(k, X_k^{t,x,*};u_k^{t,x,*}; \bar{u}_k)=\lim_{\lambda \downarrow 0}\frac{\bar{J}(k, X^{t,x,*}_k; u^{t,x,*}_k+\lambda \bar{u}_k)-\bar{J}(k, X^{t,x,*}_k; u^{t,x,*}_k)}{\lambda}\nonumber\\
&&= 2\langle Q_{k} X^{k,*}, Y^k\rangle_{\mathbb{T}_k}+2\langle R_{k,k} u_k^{t,x,*}, \bar{u}_k\rangle_k+2\langle G_{k} X_N^{k,*}, Y_N^k\rangle_{N}.
\end{eqnarray}
Similarly, the second order directional derivative with the direction $(\bar{u}_k, \hat{u}_k)$ is given by
\begin{eqnarray*}
&&d^2\bar{J}(k, X_k^{t,x,*};u_k^{t,x,*}; \bar{u}_k; \hat{u}_k)=\lim_{\beta\downarrow 0}\frac{d\bar{J}(k, X^{t,x,*}_k;u_k^{t,x,*}+\beta\hat{u}_k; \bar{u}_k)-d\bar{J}(k, X_k^{t,x,*};u_k^{t,x,*}; \bar{u}_k)}{\beta}\\
&&= 2\langle Q_{k} \hat{Y}^{k}, Y^k\rangle_{\mathbb{T}_k}+2\langle R_{k,k} \hat{u}_k, \bar{u}_k\rangle_k+2\langle G_{k} \hat{Y}_N^{k}, Y_N^k\rangle_{N},
\end{eqnarray*}
where
\begin{eqnarray*}
\left\{\begin{array}{l}
\hat{Y}^k_{\ell+1}=A_{k,\ell}\hat{Y}^k_\ell+C_{k,\ell}\hat{Y}^k_\ell w_\ell,~~~\ell\in \mathbb{T}_{k+1}, \\
\hat{Y}^k_{k+1}=B_{k,k}\bar{u}_k+D_{k,k}\bar{u}_kw_k, \\
\hat{Y}^k_k=0.
\end{array}\right.
\end{eqnarray*}
If $\hat{u}_k=\bar{u}_k$, then
\begin{eqnarray}\label{d2J}
d^2\bar{J}(k, X^{t,x,*}_k; u^{t,x,*}_k; \bar{u}_k; \bar{u}_k)=2\langle Q_{k} {Y}^{k}, Y^k\rangle_{\mathbb{T}_k}+2\langle R_{k,k} \bar{u}_k, \bar{u}_k\rangle_k+2\langle G_{k} {Y}_N^{k}, Y_N^k\rangle_{N}.
\end{eqnarray}
Note that the righthand side of (\ref{d2J}) is independent of
$u^{t,x,*}_k$. Then, for any $u_k\in
L^2_{\mathcal{F}}(k;\mathbb{R}^m)$, we have
\begin{eqnarray}\label{d2J-2}
d^2\bar{J}(k, X^{t,x,*}_k; u_k; \bar{u}_k; \bar{u}_k)=2\langle Q_{k} {Y}^{k}, Y^k\rangle_{\mathbb{T}_k}+2\langle R_{k,k} \bar{u}_k, \bar{u}_k\rangle_k+2\langle G_{k} {Y}_N^{k}, Y_N^k\rangle_{N}.
\end{eqnarray}
Furthermore, we can show that $\bar{J}(k,X_k^{t,x,*}; u_k)$ is infinitely differentiable with respect to $u_k$ in the sense that the directional derivatives of all orders exist.
%
%
By classical results on convex analysis \cite{Ekeland}, we have the following result.


\begin{lemma}\label{Lemma-convex}
The following statements are equivalent.

(i) The map $u_k\mapsto \bar{J}(k,X_k^{t,x,*}; u_k)$ is convex.

(ii) The following holds
\begin{eqnarray*}
\inf_{\bar{u}_k\in L^2_\mathcal{F}(k; \mathbb{R}^m)}\big{[}\langle Q_{k} {Y}^{k}, Y^k\rangle_{\mathbb{T}_k}+\langle R_{k,k} \bar{u}_k, \bar{u}_k\rangle_k+\langle G_{k} {Y}_N^{k}, Y_N^k\rangle_{N}\big{]}\geq 0.
\end{eqnarray*}


\end{lemma}


\textbf{Proof of Theorem \ref{Theorem-Equivalentce-open-loop}.}
(i)$\Rightarrow$(ii). Let $u^{t,x,*}$  be an open-loop equilibrium
control of Problem (LQ) for the initial pair $(t,x)$. Since
$\bar{J}(k, X^{t,x,*}_k; u_k)$ is infinitely differentiable with
respect to $u_k$ and (\ref{d2J-2}) is independent of $u_k$, the
minimum point $u^{t,x,*}_k$ of $\bar{J}(k, X^{t,x,*}_k; u_k)$ is
fully characterized via the first and second order derivatives,
namely, $d\bar{J}(k, X^{t,x,*}_k; u^{t,x,*}_k; \bar{u}_k)=0$ and
$d^2\bar{J}(k, X^{t,x,*}_k; u_k^{t,x,*}; \bar{u}_k;\bar{u}_k)\geq 0$
hold  for any $\bar{u}_k$ in $ L^2_{\mathcal{F}}(k; \mathbb{R}^m)$.
Following Lemma \ref{Lemma-convex}, (\ref{convex}) holds. The
forward S$\Delta$E of $X^{k,*}$ is clearly solvable as $Z^{k,*}$
does not appear in this S$\Delta$E. After obtaining $X^{k,*}$ and
substituting $X^{k,*}$ into the backward S$\Delta$E, we then have
$Z^{k,*}$. This means that the FBS$\Delta$E (\ref{system-adjoint})
admits a solution $(X^{k,*}, Z^{k,*})$. Furthermore, by
(\ref{first-order}), one can get
\begin{eqnarray}\label{Theorem-Equivalentce-open-loop-1}
&&\frac{1}{2}d\bar{J}(k,X_k^{t,x,*};u^{t,x,*}_k; \bar{u}_k)=\sum_{\ell=k}^{N-1}\mathbb{E}\big{[}(X_\ell^{k,*})^TQ_{k,\ell}Y_\ell^k\big{]}+\mathbb{E}\big{[}(u_k^{t,x,*})^TR_{k,k}\bar{u}_k \big{]}+\mathbb{E}\big{[}(X_N^{k,*})^TG_kY_N^t \big{]}\nonumber \\
&&=\sum_{\ell=k}^{N-1}\mathbb{E}\big{[}(X_\ell^{k,*})^TQ_{k,\ell}Y_\ell^k\big{]}+\mathbb{E}\big{[}(u_k^{t,x,*})^TR_{k,k}\bar{u}_k \big{]}+\sum_{\ell=k}^{N-1}\mathbb{E}\big{[}(Z^{k,*}_{\ell+1})^TY_{\ell+1}^k-(Z^{k,*}_{\ell})^TY_{\ell}^k \big{]}\nonumber \\
&&=\sum_{\ell=k}^{N-1}\mathbb{E}\Big{\{}\Big{[}A_{k,\ell}^T\mathbb{E}(Z_{\ell+1}^{k,*}|\mathcal{F}_{\ell-1})+Q_{k,\ell}X_\ell^{k,*}+C_{k,\ell}^T\mathbb{E}(Z_{\ell+1}^{k,*}w_\ell|\mathcal{F}_{\ell-1})-Z_\ell^{k,*}\Big{]}^TY_\ell^k\Big{\}}\nonumber \\
&&\hphantom{=}+\mathbb{E}\Big{\{}\Big{[}R_{k,k}u_k^{t,x,*}+B_{k,k}^T\mathbb{E}(Z_{k+1}^{k,*}|\mathcal{F}_{k-1})+D_{k,k}^T\mathbb{E}(Z_{k+1}^{k,*}w_t|\mathcal{F}_{k-1}) \Big{]}^T\bar{u}_k \Big{\}}.
\end{eqnarray}
As $d\bar{J}(k, X_k^{t,x,*};u_k^{t,x,*}; \bar{u}_k)=0$ for all $\bar{u}_k\in L^2_\mathcal{F}(k; \mathbb{R}^m)$ and the FBS$\Delta$E (\ref{system-adjoint}) is solvable, the stationary condition (\ref{stationary-condition}) follows.

(ii)$\Rightarrow$(i). Note that (\ref{Theorem-Equivalentce-open-loop-1}), (\ref{system-adjoint}) and (\ref{stationary-condition}), we have
\begin{eqnarray*}
d\bar{J}(k,X_k^{t,x,*};u^{t,x,*}_k; \bar{u}_k)=0,~~~\forall
\bar{u}_k\in L^2_\mathcal{F}(k; \mathbb{R}^m),
\end{eqnarray*}
which together with (\ref{convex}) implies that $u^{t,x,*}_k$ is a
minimizer of $\bar{J}(k,X_k^{t,x,*};u_k)$ over $ L^2_\mathcal{F}(k;
\mathbb{R}^m)$. Thus, (\ref{open-loop-equilibrium-4}), equivalently
(\ref{open-loop-equilibrium}), holds for $k\in \mathbb{T}_t$. This
proves the conclusion. \hfill$\square$

\subsection*{B. Proof of Theorem \ref{Theorem-Equivalentce-closed-loop}}

As (\ref{defi-closed-loop}) holds for any $u_k\in L^2_\mathcal{F}(k;
\mathbb{R}^m)$, in Definition \ref{definition-closed-loop}
(\ref{defi-closed-loop}) can be equivalently replaced by
\begin{eqnarray}\label{defi-closed-loop-1}
J\big{(}k, X_k^{t,x,*}; (\Phi X^{k,\Phi})|_{\mathbb{T}_k}\big{)}\leq J\big{(}k, X_k^{t,x,*}; (u_k+\Phi_k\overline{X}^{k,u_k,\Phi}_k,(\Phi \overline{X}^{k,u_k,\Phi})|_{\mathbb{T}_{k+1}})\big{)},
\end{eqnarray}
where
\begin{eqnarray*}
\label{equilibrium-state-closed-4}
\left\{\begin{array}{l}
\overline{X}^{k,u_k,\Phi}_{\ell+1} =\big{(}A_{k,\ell}+B_{k,\ell}\Phi_\ell\big{)} \overline{X}^{k,u_k,\Phi}_{\ell}+\big{(}C_{k,\ell}+D_{k,\ell}\Phi_\ell\big{)}\overline{X}^{k,u_k,\Phi}_\ell w_\ell,~~~\ell\in \mathbb{T}_{k+1},\\[1mm]
\overline{X}^{k,u_k,\Phi}_{k+1} =\big{(}A_{k,k}+B_{k,k}\Phi_k\big{)}\overline{X}^{k,u_k,\Phi}_k+B_{k,k}u_k+\big{[}\big{(}C_{k,k}+D_{k,k}\Phi_k\big{)}\overline{X}^{k,u_k,\Phi}_k+D_{k,k}u_k\big{]}w_k,\\[1mm]
\overline{X}^{k,u_k,\Phi}_{k} = X^{t,x,*}_k.
\end{array}\right.
\end{eqnarray*}
For any $u\in L_\mathcal{F}^2(\mathbb{T}; \mathbb{R}^m)$, $k\in
\mathbb{T}$ and $y\in L_\mathcal{F}^2(k; \mathbb{R}^m)$, let
\begin{eqnarray}\label{cost-4}
\widehat{J}(k,y; u|_{\mathbb{T}_k})= \mathbb{E}\big{[}({X}_N^k)^TG_{k}{X}^k_N\big{]}+\sum_{\ell=k}^{N-1}\mathbb{E}\big{[}({X}_\ell^k)^TQ_{k,\ell}{X}^k_{\ell}+ (\Phi_\ell{X}_\ell^k+u_{\ell})^TR_{k,\ell}(\Phi_\ell{X}_\ell^k+u_{\ell})\big{]}, 
\end{eqnarray}
where
\begin{eqnarray*}\label{system-4}
\left\{
\begin{array}{l}
{X}^k_{\ell+1} = \big{(}A_{k,\ell}+B_{k,\ell}\Phi_\ell\big{)}{X}_\ell^k+B_{k,\ell}u_{\ell}+\big{[}\big{(}C_{k,\ell}+D_{k,\ell}\Phi_\ell\big{)}{X}^k_{\ell} +D_{k,\ell}u_{\ell}\big{]}w_{\ell},\\[1mm]
X^k_k=y,~~ {\ell} \in  \mathbb{T}_k.
\end{array}
\right.
\end{eqnarray*}
Hence, (\ref{defi-closed-loop-1}) reads as
\begin{eqnarray}\label{defi-closed-loop-4}
\widehat{J}(k,X_k^{t,x,*}; 0|_{\mathbb{T}_k})\leq \widehat{J}\big{(}k, X_k^{t,x,*}; \big{(}u_k,0|_{\mathbb{T}_{k+1}}\big{)}\big{)},
\end{eqnarray}
where $0|_{\mathbb{T}_k}$ is understood as $u|_{\mathbb{T}_k}=\{u_k,\cdots, u_{N-1}\}$ with $u_\ell=0, \ell\in \mathbb{T}_k$, and similar meaning holds for $0|_{\mathbb{T}_{k+1}}$.
Therefore, if $\Phi$ is a linear feedback equilibrium strategy of Problem (LQ), then $0|_{\mathbb{T}}$ will be an open-loop equilibrium control of the time-inconsistent stochastic LQ problem corresponding to (\ref{system-1}) and (\ref{cost-4}).
The cost functional (\ref{cost-4}) is differen from (\ref{cost-1}), as in (\ref{cost-4}) crossing terms between $X_\ell^k$ and $u_\ell$ appear. The following is to mimic the proof of Theorem \ref{Theorem-Equivalentce-open-loop}.

For any given $\Phi=\{\Phi_0,\cdots, \Phi_{N-1}\}$ with $\Phi_t\in
\mathbb{R}^{m\times n}, t\in \mathbb{T}$, denote above
$\widehat{J}\big{(}k, X_k^{t,x,*};
(u_k,0|_{\mathbb{T}_{k+1}})\big{)}$ and
$\widehat{{J}}(k,X_k^{t,x,*}; 0|_{\mathbb{T}_k})$ by
$\widetilde{J}\big{(}k, X_k^{t,x,*}; u_k\big{)}$ and
$\widetilde{J}(k,X_k^{t,x,*}; 0)$, respectively. We now calculate
the  first two orders directional derivatives of
$\widetilde{J}\big{(}k, X_k^{t,x,*}; u_k\big{)}$.
Note that
\begin{eqnarray*}
&&\widetilde{J}(k,X_k^{t,x,*}; u_k)=\langle Q_k\overline{X}^{k,u_k,\Phi}, \overline{X}^{k,u_k,\Phi}\rangle_{\mathbb{T}_k}+\langle R_{k,k} (u_k+\Phi_k\overline{X}^{k,u_k,\Phi}_k), u_k+\Phi_k\overline{X}^{k,u_k,\Phi}_k\rangle_k\\[1mm]
&&\hphantom{\widetilde{J}(k,X_k^{t,x,*}; u_k)=}+\langle R_{k+}(\Phi \overline{X}^{k,u_k,\Phi})|_{\mathbb{T}_{k+1}}, (\Phi \overline{X}^{k,u_k,\Phi})|_{\mathbb{T}_k}\rangle_{\mathbb{T}_k}+\langle G_k\overline{X}^{k,u_k,\Phi}_{N}, \overline{X}^{k,u_k,\Phi}_{N}\rangle_{N},
\end{eqnarray*}
and
\begin{eqnarray*}
&&\widetilde{J}\big{(}k, X_k^{t,x,*}; u_k+\lambda \bar{u}_k\big{)}=\langle Q_k\overline{X}^{k,u_,\bar{u}_k,\Phi,\lambda}, \overline{X}^{k,u_k,\bar{u}_k,\Phi,\lambda}\rangle_{\mathbb{T}_k}\\[1mm]
&&\hphantom{\widetilde{J}\big{(}k, X_k^{t,x,*}; u_k+\lambda \bar{u}_k\big{)}=}+\langle R_{k+} (\Phi \overline{X}^{k,u_k,\bar{u}_k,\Phi,\lambda})|_{\mathbb{T}_{k+1}}, (\Phi \overline{X}^{k,u_k,\bar{u}_k,\Phi,\lambda})|_{\mathbb{T}_{k+1}}\rangle_{\mathbb{T}_{k+1}}\\[1mm]
&&\hphantom{\widetilde{J}\big{(}k, X_k^{t,x,*}; u_k+\lambda \bar{u}_k\big{)}=}+\langle R_{k,k}\big{(}u_k+\Phi_k \overline{X}^{k,u_k,\bar{u}_k,\Phi,\lambda}_k+\lambda\bar{u}_k\big{)}, u_k+\Phi_k \overline{X}^{k,u_k,\bar{u}_k,\Phi,\lambda}_k\\[1mm]
&&\hphantom{\widetilde{J}\big{(}k, X_k^{t,x,*}; u_k+\lambda \bar{u}_k\big{)}=}+\lambda\bar{u}_k\rangle_k+\langle G_k\overline{X}^{k,u_k,\bar{u}_k,\Phi,\lambda}_{N}, \overline{X}^{k,u_k,\bar{u}_k,\Phi,\lambda}_{N}\rangle_{N},
\end{eqnarray*}
where $Q_k, R_{k+}$ are similarly defined as those in (\ref{notation-1}), and
\begin{eqnarray*}
\left\{\begin{array}{l}
\overline{X}^{k,u_k,\bar{u}_k,\Phi,\lambda}_{\ell+1} =\big{(}A_{k,\ell}+B_{k,\ell}\Phi_\ell\big{)} \overline{X}^{k,u_k,\bar{u}_k,\Phi,\lambda}_{\ell}+\big{(}C_{k,\ell}+D_{k,\ell}\Phi_\ell\big{)}\overline{X}^{k,u_k,\bar{u}_k,\Phi,\lambda}_\ell w_\ell,~~\ell\in \mathbb{T}_{k+1},\\[1mm]
\overline{X}^{k,u_k,\bar{u}_k,\Phi,\lambda}_{k+1} =\big{(}A_{k,k}+B_{k,k}\Phi_k\big{)}\overline{X}^{k,u_k,\bar{u}_k,\Phi,\lambda}_k+B_{k,k}\big{(}u_k+\lambda\bar{u}_k\big{)}\\[1mm]
\hphantom{\overline{X}^{k,u_k,\bar{u}_k,\Phi,\lambda}_{k+1} =}+\big{[}\big{(}C_{k,k}+D_{k,k}\Phi_k\big{)}\overline{X}^{k,u_k,\bar{u}_k,\Phi,\lambda}_k+D_{k,k}\big{(}u_k+\lambda\bar{u}_k\big{)}\big{]}w_k,\\[1mm]
\overline{X}^{k,u_k,\bar{u}_k,\Phi,\lambda}_{k} = X^{t,x,*}_k.
\end{array}\right.
\end{eqnarray*}
As $\overline{X}^{k,u_k,\Phi}_k=\overline{X}_k^{k,u_k,\bar{u}_k,\Phi,\lambda}=X_k^{t,x,*}$, we have
\begin{eqnarray*}
\left\{\begin{array}{l}
\frac{\overline{X}_{\ell+1}^{k,u_k,\bar{u}_k,\Phi,\lambda}-\overline{X}^{k,u_k,\Phi}_{\ell+1}}{\lambda}=\big{(}A_{k,\ell}+B_{k,\ell}\Phi_\ell\big{)}\frac{\overline{X}_{\ell}^{k,u_k,\bar{u}_k,\Phi,\lambda}-\overline{X}^{k,u_k,\Phi}_{\ell}}{\lambda} \\
\hphantom{\frac{\overline{X}_{\ell+1}^{k,u_k,\bar{u}_k,\Phi,\lambda}-\overline{X}^{k,u_k,\Phi}_{\ell+1}}{\lambda}=}+\big{(}C_{k,\ell}+D_{k,\ell}\Phi_\ell\big{)}\frac{\overline{X}_{\ell}^{k,u_k,\bar{u}_k,\Phi,\lambda}-\overline{X}^{k,u_k,\Phi}_{\ell}}{\lambda}w_\ell,\\[2mm]
\frac{\overline{X}_{k+1}^{k,u_k,\bar{u}_k,\Phi,\lambda}-\overline{X}^{k,u_k,\Phi}_{k+1}}{\lambda}=B_{k,k}\bar{u}_k+D_{k,k}\bar{u}_kw_k,\\[2mm]
\frac{\overline{X}_{k}^{k,u_k,\bar{u}_k,\Phi,\lambda}-\overline{X}^{k,u_k,\Phi}_{k}}{\lambda}=0,~~~\ell\in \mathbb{T}_{k+1}.
\end{array}\right.
\end{eqnarray*}
Denote $\frac{\overline{X}_{\ell}^{k,u_k,\bar{u}_k,\Phi,\lambda}-\overline{X}^{k,\Phi}_{\ell}}{\lambda} $ by ${Y}^{k,\bar{u}_k,\Phi}_k$, which is independent of $u_k$ and $\lambda$. Then,
\begin{eqnarray}\label{system-y-closed}
\left\{\begin{array}{l}
{Y}^{k,\bar{u}_k,\Phi}_{\ell+1}=\big{(}A_{k,\ell}+B_{k,\ell}\Phi_\ell\big{)}{Y}^{k,\bar{u}_k,\Phi}_\ell+\big{(}C_{k,\ell}+D_{k,\ell}\Phi_\ell\big{)}{Y}^{k,\bar{u}_k,\Phi}_\ell w_\ell,~~~\ell\in \mathbb{T}_{k+1},\\[1mm]
Y^{k,\bar{u}_k,\Phi}_{k+1}=B_{k,k}\bar{u}_k+D_{k,k}\bar{u}_kw_k,\\[1mm]
{Y}^{k,\bar{u}_k,\Phi}_k=0.
\end{array}\right.
\end{eqnarray}
For any $\ell\in \mathbb{T}_k$, we have
\begin{eqnarray}\label{bar-X}
\overline{X}_\ell^{k,u_k,\bar{u}_k,\Phi,\lambda}=\overline{X}^{k,u_k,\Phi}_\ell+\lambda Y^{k,\bar{u}_k,\Phi}_\ell.
\end{eqnarray}
Similar to the deviation of (\ref{first-order}), we have
\begin{eqnarray}\label{first-order-closed}
&&d\widetilde{J}\big{(}k, X_k^{t,x,*}; u_k;\bar{u}_k\big{)}=\lim_{\lambda \downarrow 0}\frac{\widetilde{J}\big{(}k, X_k^{t,x,*}; u_k+\lambda \bar{u}_k\big{)}-\widetilde{J}\big{(}k, X_k^{t,x,*}; u_k\big{)}}{\lambda}\nonumber\\
&&= 2\langle Q_{k} \overline{X}^{k,u_k,\Phi}, Y^{k,\bar{u}_k,\Phi}\rangle_{\mathbb{T}_k}+2\langle R_{k,k} (u_k+\Phi_k \overline{X}_k^{k,u_k,\Phi}), \bar{u}_k\rangle_k\nonumber\\[1mm]
&&\hphantom{=}+2\langle R_{k+}(\Phi \overline{X}^{k,u_k,\Phi})|_{\mathbb{T}_{k+1}},  (\Phi Y^{k,\bar{u}_k,\Phi})|_{\mathbb{T}_{k+1}} \rangle_{\mathbb{T}_{k+1}}\nonumber\\
&&\hphantom{=}+2\langle G_{k} \overline{X}_N^{k,u_k,\Phi}, Y_N^{k,\bar{u}_k,\Phi}\rangle_{N},
\end{eqnarray}
and for any $\hat{{u}}_k\in L^2_\mathcal{F}(k; \mathbb{R}^m)$
\begin{eqnarray*}\label{second-order-2}
&&d^2\widetilde{J}\big{(}k, X_k^{t,x,*}; u_k; \bar{u}_k;\hat{u}_k\big{)}=\lim_{\lambda \downarrow 0}\frac{\widetilde{J}\big{(}k, X_k^{t,x,*}; u_k+\beta \hat{u}_k;\bar{u}_k\big{)}-\widetilde{J}\big{(}k, X_k^{t,x,*}; u_k;\bar{u}_k\big{)}}{\beta}\nonumber\\[1mm]
&&= 2\langle Q_k \hat{Y}^{k,\hat{u}_k,\Phi}, Y^{k,\bar{u}_k, \Phi}\rangle_{\mathbb{T}_k}+2\langle R_{k,k}\hat{u}_k, \bar{u}_k\rangle_k+2\langle R_{k+}(\Phi \hat{Y}^{k,\hat{u}_k,\Phi})|_{\mathbb{T}_{k+1}}, (\Phi Y^{k,\bar{u}_k, \Phi})|_{\mathbb{T}_{k+1}}\rangle_{\mathbb{T}_{k+1}}\\[1mm]
&&\hphantom{=}+2\langle G_k \hat{Y}_N^{k,\hat{u}_k, \Phi}, Y_N^{k,\bar{u}_k, \Phi}\rangle_N,
\end{eqnarray*}
where
\begin{eqnarray*}\label{system-hat-y}
\left\{
\begin{array}{l}
\hat{Y}^{k,\hat{u}_k,\Phi}_{\ell+1}=\big{(}A_{k,\ell}+B_{k,\ell}\Phi_\ell\big{)}\hat{Y}^{k,\hat{u}_k,\Phi}_\ell+\big{(}C_{k,\ell}+D_{k,\ell}\Phi_\ell\big{)}\hat{Y}^{k,\hat{u}_k,\Phi}_\ell w_\ell,~~~\ell\in \mathbb{T}_{k+1},\\[1mm]
Y^{k,\hat{u}_k,\Phi}_{k+1}=B_{k,k}\hat{u}_k+D_{k,k}\hat{u}_kw_k,\\[1mm]
\hat{Y}^{k,\hat{u}_k,\Phi}_k=0.
\end{array}
\right.
\end{eqnarray*}
If $\hat{u}_k=\bar{u}_k$, then we have
\begin{eqnarray}\label{d2J-closed}
&&d^2\widetilde{J}\big{(}k, X_k^{t,x,*}; u_k;\bar{u}_k; \bar{u}_k\big{)}=2\langle Q_k {Y}^{k,\bar{u}_k,\Phi}, Y^{k,\bar{u}_k, \Phi}\rangle_{\mathbb{T}_k}+2\langle R_{k,k}\bar{u}_k, \bar{u}_k\rangle_k\nonumber \\[1mm]
&&\hphantom{d^2\widetilde{J}\big{(}k, X_k^{t,x,*}; u_k;\bar{u}_k; \bar{u}_k\big{)}=}+2\langle R_{k+}(\Phi {Y}^{k,\bar{u}_k,\Phi})|_{\mathbb{T}_{k+1}}, (\Phi Y^{k,\bar{u}_k, \Phi})|_{\mathbb{T}_{k+1}}\rangle_{\mathbb{T}_{k+1}}\nonumber \\[1mm]
&&\hphantom{d^2\widetilde{J}\big{(}k, X_k^{t,x,*}; u_k;\bar{u}_k; \bar{u}_k\big{)}=}+2\langle G_k {Y}_N^{k,\bar{u}_k, \Phi}, Y_N^{k,\bar{u}_k, \Phi}\rangle_N,
\end{eqnarray}
which is independent of $u_k$.

\textbf{{Proof of Theorem \ref{Theorem-Equivalentce-closed-loop}}}.
(i)$\Rightarrow$(ii). Let $\Phi$  be a linear feedback equilibrium strategy of Problem (LQ).
As noted above, $0|_{\mathbb{T}}$ is an open-loop equilibrium control of the stochastic LQ problem corresponding to (\ref{system-1}) and (\ref{cost-4}). Hence, for all $\bar{u}_k \in L_{\mathcal{F}}^2(k; \mathbb{R}^m)$ we have
\begin{eqnarray}\label{defi-closed-loop-3}
d\widetilde{J}\big{(}k, X_k^{t,x,*}; 0;\bar{u}_k\big{)}=0,~~~
d^2\widetilde{J}\big{(}k, X_k^{t,x,*}; 0; \bar{u}_k;\bar{u}_k\big{)}\geq 0.
\end{eqnarray}
Note that $\overline{X}^{k,u_k,\Phi}_\ell=X^{k,\Phi}_\ell, \ell\in \mathbb{T}_k$ if $u_k=0$. Hence, it holds that
\begin{eqnarray}\label{first-order-closed-2}
&&0=\frac{1}{2}d{\widetilde{J}}(k, X_k^{t,x,*}; 0; \bar{u}_k)\nonumber \\ [1mm]
&&\hphantom{0}=\sum_{\ell=k}^{N-1}\mathbb{E}\Big{[}\Big{(}\big{(}A_{k,\ell}+B_{k,\ell}\Phi_\ell\big{)}^T\mathbb{E}(Z_{\ell+1}^{k,\Phi}|\mathcal{F}_{\ell-1})+\big{(}C_{k,\ell}+D_{k,\ell}\Phi_\ell\big{)}^T\mathbb{E}(Z_{\ell+1}^{k,\Phi}w_\ell|\mathcal{F}_{\ell-1})\nonumber \\[1mm]
&&\hphantom{0=}+\big{(}\Phi_\ell^TR_{k,\ell}\Phi_\ell+Q_{k,\ell}\big{)}X_\ell^{k,\Phi}-Z_\ell^{k,\Phi}\Big{)}^TY_\ell^{k,\bar{u}_k,\Phi}\Big{]} +\mathbb{E}\Big{[}\Big{(}R_{k,k}\Phi_k X_k^{k,\Phi}+B_{k,k}^T\mathbb{E}(Z_{k+1}^{k,\Phi}|\mathcal{F}_{k-1})\nonumber  \\[1mm]
&&\hphantom{0=}+D_{k,k}^T\mathbb{E}(Z_{k+1}^{k,\Phi}w_k|\mathcal{F}_{k-1}) \Big{)}^T\bar{u}_k \Big{]}\nonumber \\[1mm]
&&\hphantom{0}=\mathbb{E}\Big{[}\Big{(}R_{k,k}\Phi_k X_k^{k,\Phi}+B_{k,k}^T\mathbb{E}(Z_{k+1}^{k,\Phi}|\mathcal{F}_{k-1})+D_{k,k}^T\mathbb{E}(Z_{k+1}^{k,\Phi}w_k|\mathcal{F}_{k-1}) \Big{)}^T\bar{u}_k \Big{]}.
\end{eqnarray}
As $\bar{u}_k\in L_{\mathcal{F}}^2(k; \mathbb{R}^m)$ can be taken arbitrarily, we can get (\ref{stationary-condition-closed}). (\ref{convex-closed}) follows from (\ref{d2J-closed}) and (\ref{defi-closed-loop-3}).

(ii)$\Rightarrow$(i). Reversing the procedure of the above proof, we can achieve the conclusion. The proof is also similar to that of Theorem \ref{Theorem-Equivalentce-open-loop}. \hfill $\square$




\end{document}